\documentclass[11pt, reqno]{amsart}

\usepackage{amsmath}

\usepackage{tabu}
\usepackage{tikz}
\usepackage[margin=1.2in]{geometry}
\usepackage[matrix,arrow,curve,frame]{xy}
 \usepackage{hyperref}
\usepackage{amsthm}
\usepackage{amsfonts}
\usepackage{amssymb}
\usepackage{wasysym}
\usepackage{mathrsfs}
\usepackage{mathtools}
\usepackage{graphicx}

\definecolor{my-linkcolor}{rgb}{0.75,0,0}
\definecolor{my-citecolor}{rgb}{0.1,0.57,0}
\definecolor{my-urlcolor}{rgb}{0,0,0.75}
\hypersetup{
	colorlinks,
	linkcolor={my-linkcolor},
	citecolor={my-citecolor},
	urlcolor={my-urlcolor}
}

\title[Rigidity for braided r-categories]{Rigidity for pivotal braided r-categories of moderate growth and applications to vertex operator algebras}
 \author{Robert McRae}
\date{}

 \address{Yau Mathematical Sciences Center, Tsinghua University, Beijing 100084, China}
  \email{rhmcrae@tsinghua.edu.cn}


\newtheorem{thm}{Theorem}[section]
\newtheorem{cor}[thm]{Corollary}
\newtheorem{lem}[thm]{Lemma}
\newtheorem{prop}[thm]{Proposition}
\theoremstyle{definition}
\theoremstyle{definition}\newtheorem{rem}[thm]{Remark}
\theoremstyle{definition}\newtheorem{exam}[thm]{Example}
\theoremstyle{definition}

\newcommand{\cA}{\mathcal{A}}

\newcommand{\cC}{\mathcal{C}}

\newcommand{\til}{\widetilde}

\newcommand{\CC}{\mathbb{C}}
\newcommand{\ZZ}{\mathbb{Z}}
\newcommand{\NN}{\mathbb{N}}

\newcommand{\QQ}{\mathbb{Q}}

\newcommand{\Id}{\mathrm{Id}}

\newcommand{\tens}{\otimes}
\newcommand{\vac}{\mathbf{1}}

 \DeclareMathOperator{\im}{Im}
 \DeclareMathOperator{\coker}{Coker}

 \DeclareMathOperator{\rep}{Rep}
 \let\ker\relax
 \let\hom\relax
 \DeclareMathOperator{\ker}{Ker}
 \DeclareMathOperator{\hom}{Hom}

\begin{document}
\bibliographystyle{alpha}

\numberwithin{equation}{section}

 \begin{abstract}
Let $\mathcal{C}$ be an abelian pivotal braided r-category with simple unit object $\mathbf{1}$ over any field.  We present a ChatGPT-generated proof, prompted by Bin Gui, that any object $X$ in $\mathcal{C}$ of moderate growth, in the sense that $\dim\mathrm{End}(X^{\otimes n})<n!$ for some $n\in\mathbb{Z}_{\geq 1}$, is rigid. The proof crucially uses ideas from the recent proof of Etingof and Penneys that $X$ is rigid if it is also non-negligible, that is, $X\otimes Y$ contains $\mathbf{1}$ as a direct summand for some object $Y$. In the present paper, the non-negligible assumption is removed using a rigidity criterion developed previously by the author for use in non-semisimple categories. We also present a variant result (found without further AI assistance) that $X$ is also rigid if $\dim\mathrm{Hom}((DX)^{\otimes n}\otimes X^{\otimes n},\vac)<n!$ for some $n\in\mathbb{Z}_{\geq 1}$, where $D$ is the Grothendieck-Verdier duality functor of $\mathcal{C}$. Using this variant, we show that if $V$ is a simple self-contragredient vertex operator algebra such that the category $\mathcal{C}^1_V$ of $C_1$-cofinite grading-restricted generalized $V$-modules is closed under contragredient modules, then $\mathcal{C}^1_V$ is a locally finite $\CC$-linear abelian rigid braided ribbon category. We also show that if $V$ is a simple $\mathbb{Z}_{\geq 0}$-graded self-contragredient $C_2$-cofinite vertex operator algebra, then the category of all grading-restricted generalized $V$-modules is a (not necessarily semisimple) modular tensor category.
 \end{abstract}

\maketitle

\tableofcontents

\section{Introduction}

Let $\Bbbk$ be a field, and let $\cC$ be a $\Bbbk$-linear abelian monoidal category with unit object $\vac$ and tensor product bifunctor $\otimes: \cC\times\cC\rightarrow\cC$. Then $\cC$ is rigid if every object has duals in a strong sense. Specifically, for $X\in\mathrm{Ob}(\cC)$, a left dual of $X$ is an object $X^*$ with evaluation and coevaluation morphisms $e_X: X^*\otimes X\rightarrow\vac$ and $i_X:\vac\rightarrow X\otimes X^*$ such that the compositions
\begin{equation*}
X\xrightarrow{\sim} \vac\otimes X\xrightarrow{i_X\otimes\Id} (X\otimes X^*)\otimes X\xrightarrow{\sim} X\otimes (X^*\otimes X)\xrightarrow{\Id\otimes e_X} X\otimes\vac\xrightarrow{\sim} X
\end{equation*}
and
\begin{equation*}
X^*\xrightarrow{\sim} X^*\otimes \vac\xrightarrow{\Id\otimes i_X} X^*\otimes(X\otimes X^*)\xrightarrow{\sim} (X^*\otimes X)\otimes X^*\xrightarrow{e_X\otimes\Id} \vac\otimes X^*\xrightarrow{\sim} X^*
\end{equation*}
are both identities. Right duals are defined analogously, and $\cC$ is rigid if every object has a left and right dual. Many monoidal categories are obviously rigid. For example, Hopf algebras are essentially defined so that finite-dimensional representations of a Hopf algebra over $\Bbbk$ form a rigid $\Bbbk$-linear abelian monoidal category. Consequently, the categories of finite-dimensional representations of groups, Lie algebras, and quantum groups are all rigid. For other categories, existence of duals is more interesting. For example, in monoidal categories of endofunctors, where morphisms are natural transformations and the tensor product is functor composition, left and right duals are (if they exist) left and right adjoint functors.

In this paper, we focus on braided monoidal categories $\cC$, where there are natural commutativity (or braiding) isomorphisms $c_{X,Y}: X\otimes Y\rightarrow Y\otimes X$ for $X,Y\in\mathrm{Ob}(\cC)$, and where left and right duals are isomorphic as objects of $\cC$ (and both exist if one does). We do not assume $\cC$ is symmetric, that is, we do not assume $c_{Y,X} = c_{X,Y}^{-1}$ in general, so that $c_{Y,X}\circ c_{X,Y}$ may be a non-trivial automorphism of $X\otimes Y$. The braiding isomorphisms and rigidity are both essential to the well-known applications of braided monoidal categories to knot and link invariants, $3$-manifold invariants, and topological quantum field theories. Thus it is natural to ask whether there are general conditions under which a braided monoidal category is rigid.

\subsection{Main results}

In this paper, we present a ChatGPT-generated proof of rigidity for a large class of braided monoidal categories. The starting point is a weaker duality structure on a $\Bbbk$-linear abelian braided monoidal category $\cC$, namely, we assume that $\cC$ is an \textit{r-category} in the sense of \cite{BD} (or, $\cC$ is \textit{weakly rigid} in the sense of \cite[Definition A.4]{KL4}). This means in particular that there is an anti-equivalence $D:\cC\rightarrow\cC$, which we call the duality functor, such that for any $X\in\mathrm{Ob}(\cC)$, there is a morphism $e_X: DX\otimes X\rightarrow\vac$ that acts like the evaluation morphism of a left dual. This structure is weaker than rigidity, since while any rigid category is an r-category with $D$ given by left duals, in a general r-category, there is no assumption of coevaluation-like morphisms $i_X: \vac\rightarrow X\otimes DX$. An r-category is called pivotal if there is a natural isomorphism $\delta:\Id_\cC\rightarrow D^2$ satisfying suitable properties. Then the main result is Theorem \ref{thm:main-theorem} below:

\begin{thm}\label{thm:intro-main-theorem}
Let $\cC$ be a $\Bbbk$-linear abelian pivotal braided r-category with a simple unit object $\vac$, and let $X\in\mathrm{Ob}(\cC)$. If there exists $n\in\ZZ_{\geq 1}$ such that either
\begin{equation*}
\dim\mathrm{Hom}((DX)^{\otimes n}\otimes X^{\otimes n},\vac)<n!
\end{equation*}
or
\begin{equation*}
\dim\mathrm{End}(X^{\otimes n})<n!,
\end{equation*}
then $X$ is rigid with left dual $DX$.
\end{thm}

The two dimension conditions in this theorem are not obviously equivalent unless we already know that $DX$ is a left dual of $X$. But either one is straightforward to check if $\cC$ is locally finite abelian and has finitely many isomorphism classes of simple objects, so one important corollary is Theorem \ref{thm:fin-abelian-rigid} below:
\begin{cor}\label{cor:intro-finite}
Let $\cC$ be a locally finite $\Bbbk$-linear abelian pivotal braided r-category with simple unit object $\vac$ and finitely many isomorphism classes of simple objects. Then $\cC$ is rigid with left duals given by the duality functor $D$.
\end{cor}

The main idea in the proof of Theorem \ref{thm:intro-main-theorem} comes from the recent major breakthrough \cite{EP} by Etingof and Penneys on rigidity for general braided monoidal categories, where the condition $\dim\mathrm{End}(X^{\otimes n})<n!$ for some $n\in\ZZ_{\geq 1}$ is called \textit{moderate growth}. In \cite[Theorem 1.1]{EP}, Etingof and Penneys showed that if $\cC$ is closed under direct sums and summands and has finite-dimensional morphism spaces, and if $\mathrm{End}(\vac)=\Bbbk$, then $X\in\mathrm{Ob}(\cC)$ is rigid if it has moderate growth and is also \textit{non-negligible}. This last condition means that $X$ is indecomposable and there exists $Y\in\mathrm{Ob}(\cC)$ (which we may assume is also indecomposable) such that $X\otimes Y$ contains $\vac$ as a direct summand.

The proof of \cite[Theorem 1.1]{EP} uses the observation that, through the braiding $c$, the symmetric group $S_n$ maps into $\mathrm{End}(X^{\otimes n})$ as a subset for any $X\in\mathrm{Ob}(\cC)$ and $n\in\ZZ_{\geq 1}$. Assuming $X$ is non-negligible, one has morphisms $\til{e}_X: X\otimes Y\rightarrow\vac$ and $i_X:\vac\rightarrow X\otimes Y$ for some $Y\in\mathrm{Ob}(\cC)$ such that $\til{e}_X\circ i_X=\Id$, and using these maps, one defines a 
bilinear form 
\begin{equation*}
\phi_n: \mathrm{End}(X^{\otimes n})\times\mathrm{End}(X^{\otimes n})\longrightarrow\mathrm{End}(\vac)=\Bbbk
\end{equation*}
such that $\phi_n(f,g)$ is the composition
\begin{equation*}
\vac\longrightarrow X^{\otimes n}\otimes Y^{\otimes n} \xrightarrow{(f\circ g)\otimes\Id} X^{\otimes n}\otimes Y^{\otimes n}\longrightarrow\vac.
\end{equation*}
Etingof and Penneys then show that if $\til{e}_X$ and $i_X$ \textit{cannot} be modified into an evaluation and coevaluation for $X$ (with left dual $Y$) then the morphisms $f_\sigma\in\mathrm{End}(X^{\otimes n})$ coming from permutations $\sigma\in S_n$ satisfy
\begin{equation*}
\phi_n(f_\sigma^{-1},f_\tau)=\delta_{\sigma,\tau},
\end{equation*}
and thus the $f_\sigma$ are linearly independent. Consequently, if a non-negligible $X\in\mathrm{Ob}(\cC)$ is not rigid, then $X$ cannot have moderate growth.

While \cite[Theorem 1.1]{EP} is already extremely useful, especially for semisimple braided r-categories where every non-zero object is non-negligible, the non-negligibility assumption can be a serious obstacle in non-semisimple categories. Thus Theorem \ref{thm:intro-main-theorem} now removes this obstacle at the cost of the relatively mild assumption that $\cC$ is a pivotal r-category. However, this improvement of \cite[Theorem 1.1]{EP} is not completely straightforward. Indeed, if $X$ is an object of a pivotal braided r-category, then it is natural to take the object $Y$ in \cite{EP} to be $DX$, since there is a natural map
\begin{equation*}
\til{e}_X = e_{DX}\circ(\delta_X\otimes\Id_{DX}): X\otimes DX\longrightarrow\vac,
\end{equation*}
where $\delta_X: X\rightarrow D^2X$ is the pivotal structure; equivalently,
\begin{equation*}
\til{e}_X= e_X\circ c_{X,DX}\circ(\theta_X\otimes\Id_{DX})
\end{equation*}
where $\theta_X$ is the natural twist automorphism of $X$, which is equivalent to the pivotal structure.
 But if $X$ is negligible, then no suitable morphism $i_X: \vac\rightarrow X\otimes DX$ is guaranteed, so we cannot define a bilinear form on $\mathrm{End}(X^{\otimes n})$ as in \cite{EP}.

The solution to this problem comes from a criterion for rigidity in possibly non-semisimple pivotal braided r-categories proved in \cite{McR}. To start, one shows that the pairing
\begin{equation*}
\til{e}_X\circ(\Id_X\otimes e_X\otimes \Id_{DX}): (X\otimes DX)\otimes (X\otimes DX)\longrightarrow\vac
\end{equation*}
induces a natural morphism
\begin{equation*}
\Phi_X: X\otimes DX\longrightarrow D(X\otimes DX)
\end{equation*}
such that $X$ is rigid if and only if $\Phi_X$ is an isomorphism (see for example \cite[Corollary 3.6]{BD}), in which case $DX$ is the left dual of $X$, the evaluation is $e_X$, and the coevaluation $i_X$ is the composition
\begin{equation*}
\vac\xrightarrow{\sim} D\vac\xrightarrow{D\til{e}_X} D(X\otimes DX)\xrightarrow{\Phi_X^{-1}} X\otimes DX.
\end{equation*}
In general, it may be no easier to construct $\Phi_X^{-1}$ explicitly than it is to construct a coevaluation candidate $i_X$. However, as proved in Theorem \ref{thm:rigidity-criterion} below, it follows straightforwardly from Theorem 3.4 and Proposition 3.5 of \cite{McR} that $\Phi_X$ is an isomorphism, and thus $X$ is rigid, if $\ker\Phi_X\subseteq\ker\til{e}_X$. That is, letting $k$ denote the inclusion of $\ker\Phi_X$ into $X\otimes DX$, if $X$ is \textit{not} rigid, then $\til{e}_X\circ k\neq 0$.

Using this rigidity criterion, the moderate growth methods of \cite{EP} generalize to prove rigidity of possibly negligible objects in pivotal braided r-categories $\cC$. Indeed, suppose $X\in\mathrm{Ob}(\cC)$ is not rigid, so that $\til{e}_X\circ k\neq 0$, where $k$ is the kernel morphism of $\Phi_X$. We can then define a bilinear map
\begin{equation*}
\phi_n: \mathrm{End}(X^{\otimes n})\times\mathrm{End}(X^{\otimes n})\longrightarrow\mathrm{Hom}(K^{\otimes n},\vac),
\end{equation*}
where $K=\ker\Phi_X$,
such that $\phi_n(f,g)$ is the composition
\begin{equation*}
K^{\otimes n} \xrightarrow{k^{\otimes n}} (X\otimes DX)^{\otimes n} \xrightarrow{\sim} (DX)^{\otimes n}\otimes X^{\otimes n}\xrightarrow{\Id\otimes (f\circ g)} (DX)^{\otimes n}\otimes X^{\otimes n}\longrightarrow\vac,
\end{equation*}
where the second arrow is a suitable composition of braiding and twist automorphisms and the last arrow is obtained using $e_X$. One can now map $\sigma\in S_n$ to the automorphism $f_\sigma\in\mathrm{End}(X^{\otimes n})$ as in \cite{EP} and show that
\begin{equation*}
\phi_n(f_\sigma^{-1},f_\tau) =\delta_{\sigma,\tau} (\til{e}_X\circ k)^{\otimes n}
\end{equation*}
for $\sigma,\tau\in S_n$. Assuming that $\vac$ is simple and $X$ is not rigid, the morphism $\til{e}_X\circ k$ is surjective, so that $(\til{e}_X\circ k)^{\otimes n}$ is also surjective and non-zero by right exactness of $\otimes$. It then follows as in \cite{EP} that the morphisms $\lbrace f_\sigma\mid \sigma\in S_n\rbrace $ are linearly independent and thus $X$ does not have moderate growth. One shows similarly that the morphisms
\begin{equation*}
(DX)^{\otimes n}\otimes X^{\otimes n} \xrightarrow{\Id\otimes f_\sigma} (DX)^{\otimes n}\otimes X^{\otimes n} \longrightarrow\vac
\end{equation*}
are linearly independent if $X$ is not rigid, so $\dim\hom((DX)^{\otimes n}\otimes X^{\otimes n},\vac)\geq n!$ as well. This proves Theorem \ref{thm:intro-main-theorem}.

\subsection{Application to vertex operator algebras}

The original motivation for Theorem \ref{thm:intro-main-theorem} and Corollary \ref{cor:intro-finite} is the problem of rigidity for braided monoidal categories of modules for vertex operator algebras, which are algebraic structures arising from the mathematics of two-dimensional conformal quantum field theory in physics, where they appear as the chiral algebras of conformal field theories. The physics notion of fusion of two representations of a chiral algebra motivates the definition of tensor product $\boxtimes$ of modules for a vertex operator algebra, and the associativity of $\boxtimes$ amounts to the existence of operator product expansions of quantum fields in a conformal field theory. The work of Huang, Lepowsky, and Zhang in \cite{HLZ1}-\cite{HLZ8} gives sufficient conditions under which $\boxtimes$ gives categories of modules for a vertex operator algebra braided monoidal structure, and by checking these sufficient conditions, Huang has shown recently in \cite{Hu} that the category of so-called ``$C_1$-cofinite'' modules for any vertex operator algebra is naturally a braided monoidal category. Thus vertex operator algebras are a rich source of braided monoidal categories, and hence also of knot and link invariants, $3$-manifold invariants, and topological quantum field theories if these categories are rigid.

However, it is not usually obvious whether a braided monoidal category of modules for a vertex operator algebra is rigid, since the rigidity axioms involve associativity isomorphisms, and for vertex operator algebras, the associativity of $\boxtimes$ is difficult to describe explicitly. The natural candidate for the dual of a module $W$ for a vertex operator algebra $V$ is its contragredient $W'$, which is a natural module structure on the graded dual vector space of $W$ \cite{FHL}. If a braided monoidal category $\cC$ of $V$-modules is closed under contragredients, then it is easy to see that $\cC$ is a ribbon Grothendieck-Verdier category in the sense of \cite{BD} (see for example \cite[Theorem 2.12]{ALSW}). In particular, if $V\cong V'$ as a $V$-module, then $\cC$ is a pivotal braided r-category, and we can hope that $\cC$ is rigid.

Probably the first rigidity result for braided monoidal categories coming from vertex operator algebras was Kazhdan and Lusztig's proof in \cite{KL4} of rigidity for certain categories of modules for affine Lie algebras, which turn out to be the categories of $C_1$-cofinite modules for suitable affine vertex operator algebras.
The first general rigidity result for vertex operator algebras was then Huang's theorem \cite{Hu-Verlinde, Hu-rig-mod} that the representation category of a so-called ``strongly rational'' vertex operator algebra is a rigid, and indeed modular, braided tensor category. A vertex operator algebra is strongly rational if it satisfies the strong finiteness condition of $C_2$-cofiniteness introduced by Zhu \cite{Zh}, among other relatively minor conditions, and if its representation category is semisimple; such vertex operator algebras arise in rational conformal field theory in physics. The rigidity proofs of Kazhdan-Lusztig and Huang are complicated but can be simplified using \cite{EP} (see \cite[Section 4.3]{EP} and the discussion following \cite[Corollary 1.3]{EP}).

Dropping the semisimplicity assumption from the definition of strongly rational vertex operator algebra, we get simple $\ZZ_{\geq 0}$-graded self-contragredient $C_2$-cofinite vertex operator algebras. It was conjectured in \cite{Hu-conj} that the representation category $\rep(V)$ of such a vertex operator algebra $V$ should be a (possibly non-semisimple) modular tensor category, and in particular should be rigid. It was then proved in \cite{McR} that indeed $\rep(V)$ is modular if rigidity holds. Before now, rigidity has been checked for most well-known examples, including the triplet $W$-algebras \cite{TW}, their finite cyclic orbifolds \cite{CMY4}, and the symplectic fermion vertex operator algebras \cite{McR-Deligne}. But now, Huang's conjecture in all cases follows from Corollary \ref{cor:intro-finite} (see Theorem \ref{thm:C2-cofinite-rigidity} below):
\begin{thm}
If $V$ is a simple $\ZZ_{\geq 0}$-graded self-contragredient $C_2$-cofinite vertex operator algebra, then the category $\rep(V)$ of grading-restricted generalized $V$-modules is a (not necessarily semisimple) modular tensor category.
\end{thm}

Most vertex operator algebras are not $C_2$-cofinite and have infinitely many simple modules, so Corollary \ref{cor:intro-finite} does not apply. But rigidity is still expected for braided monoidal categories from non-$C_2$-cofinite vertex operator algebras under reasonable conditions. Rigidity has already been proved in many examples, including in categories coming from many affine vertex operator algebras \cite{KL4, CY, MY, CMY6, NORW}, the affine vertex operator superalgebra of $\mathfrak{gl}_{1\vert 1}$ \cite{CMY3}, many Virasoro vertex operator algebras \cite{CJORY, MY1, MY2}, and the singlet $W$-algebras \cite{CMY2,CMY4}. Rigidity has also been partially proved for $C_1$-cofinite modules of the $N=1$ super Virasoro vertex operator superalgebras \cite{CMY1,CMOY} and affine $W$-algebras associated to simply-laced Lie algebras at irrational levels \cite{CDN}. But before now, there has been no general rigidity theorem for non-$C_2$-cofinite vertex operator algebras. In Theorem \ref{thm:C1-cofinite-rigidity} below, we show that if $W$ is a $C_1$-cofinite module for any simple self-contragredient vertex operator algebra $V$ such that the category of $C_1$-cofinite $V$-module is closed under contragredients, then there exist constants $A,B>0$ such that
\begin{equation*}
\dim\hom((W')^{\boxtimes n}\boxtimes W^{\boxtimes n}, V)\leq A\cdot B^n,
\end{equation*}
for all $n\in\ZZ_{\geq 1}$, and consequently Theorem \ref{thm:intro-main-theorem} implies:
\begin{thm}
Let $V$ be a simple $\ZZ$-graded vertex operator algebra such that $V\cong V'$ as a $V$-module, and assume that the category $\cC_V^1$ of $C_1$-cofinite grading-restricted generalized $V$-modules is closed under contragredients. Then $\cC_V^1$ is a locally finite $\CC$-linear abelian braided ribbon category. In particular, $\cC_V^1$ is rigid with duals given by contragredient modules.
\end{thm}

This theorem also applies to $C_1$-cofinite modules for suitable vertex operator superalgebras. As an application, we get rigidity of $C_1$-cofinite modules for all simple Virasoro and all simple $N=1$ super Virasoro vertex operator (super)algebras. We also get rigidity of the Kazhdan-Lusztig category $KL^k(\mathfrak{g})$ of level-$k$ modules for the affine Lie algebra $\widehat{\mathfrak{g}}$, defined in \cite[Definition 2.15]{KL1}, for any simple Lie algebra $\mathfrak{g}$ and any level $k$ such that $k+h^\vee\notin\QQ_{\geq 0}$, where $h^\vee$ is the dual Coxeter number of $\mathfrak{g}$. For such $k$, $KL^k(\mathfrak{g})$ is precisely the category of $C_1$-cofinite modules for the simple affine vertex operator algebra of $\mathfrak{g}$ at level $k$. In particular, unlike in \cite{KL1}-\cite{KL4}, we do not need to assume $\mathfrak{g}$ is simply laced, and if $k+h^\vee\in\QQ_{<0}$, we also do not need to assume, as in \cite{KL4} and \cite[Section 4.3]{EP}, that the numerator of $\vert k+h^\vee\vert$ is sufficiently large. If $k+h^\vee\in\QQ_{>0}$, then the category of $C_1$-cofinite modules for the simple affine vertex operator algebra associated to $\mathfrak{g}$ at level $k$ is still a braided monoidal category by \cite{Hu}; we will show that it is rigid and semisimple in \cite{McR3}.

\subsection{AI use disclosure and methodology}

The proof of rigidity for objects of moderate growth in an abelian pivotal braided r-category with simple unit object was generated by ChatGPT 6 Pro through prompting by the author's Tsinghua University colleague Bin Gui. First, ChatGPT was asked to prove rigidity of the representation category of a simple $\ZZ_{\geq 0}$-graded self-contragredient $C_2$-cofinite vertex operator algebra, and it produced the correct idea for the proof. Then Bin Gui asked ChatGPT to prove the result for more general r-categories, and it generated the statements and proofs of the following results:

\begin{thm}\label{thm:intro-ChatGPT}
Let $\Bbbk$ be an algebraically closed field, and let $\cC$ be a locally finite $\Bbbk$-linear abelian pivotal braided r-category with simple tensor unit. Assume that for every simple object $S$, there exists an integer $n\geq 2$ such that
\begin{equation*}
\dim\mathrm{End}(S^{\otimes n})<n!.
\end{equation*}
Then $\cC$ is rigid, and the r-category dual $DX$ can be taken as the left dual of any object $X$.
\end{thm}

\begin{cor}\label{cor:intro-ChatGPT}
Let $\Bbbk$ be an algebraically closed
field. Every locally finite $\Bbbk$-linear abelian braided pivotal r-category with simple unit
and finitely many isomorphism classes of simple objects is rigid.
\end{cor}

The ChatGPT output was then shared with the author, who checked the correctness of the proof, in particular observing ChatGPT's successful combination of the moderate growth ideas and methods from \cite{EP} with the rigidity criterion for non-semisimple pivotal r-categories in \cite{McR}. The author then undertook to write an exposition of the proof, resulting in Theorem \ref{thm:intro-main-theorem} and Corollary \ref{cor:intro-finite}, and in particular removing the assumption that $\Bbbk$ is algebraically closed. Bin Gui declined to participate in the write-up. AI tools were not used in the writing of this paper. Thus the proof of Theorem \ref{thm:intro-main-theorem} differs in some details and structure from the ChatGPT proof of Theorem \ref{thm:intro-ChatGPT} since it was not written by simply modifying the original ChatGPT output. However, the proof is essentially the same as it uses the same key ideas and methods.

Most of Section \ref{sec:rig-crit} reviews definitions and results from \cite{BD} and \cite{McR} that are used in the proof of Theorem \ref{thm:intro-main-theorem}; the observation that Theorem \ref{thm:rigidity-criterion} and Corollary \ref{cor:non-rigidity-consequence} follow from Theorem 3.4 and Proposition 3.5 of \cite{McR} is due to ChatGPT. Section \ref{sec:rigidity-proof} contains the author's exposition of the variants Theorem \ref{thm:intro-main-theorem} and Corollary \ref{cor:intro-finite} of the ChatGPT-generated Theorem \ref{thm:intro-ChatGPT} and Corollary \ref{cor:intro-ChatGPT}. In the process of rewriting the ChatGPT-generated proof, the author discovered (without further AI assistance) the variant sufficient condition
\begin{equation*}
\dim\hom((DX)^{\otimes n}\otimes X^{\otimes n},\vac)<n!
\end{equation*}
for rigidity. This variant is essential for Theorem \ref{thm:C1-cofinite-rigidity} in Section \ref{sec:VOA-app} on rigidity of $C_1$-cofinite modules for vertex operator algebras, which was proved without AI assistance. Theorem \ref{thm:C2-cofinite-rigidity} for $C_2$-cofinite vertex operator algebras was observed by ChatGPT based on Corollary \ref{cor:intro-ChatGPT}.

 The author takes responsibility for this manuscript, including for correctness of the references and of the results presented here.

\subsection*{Acknowledgments}

I would like to thank Bin Gui for sharing the ChatGPT-generated proofs of Theorem \ref{thm:intro-ChatGPT} and Corollary \ref{cor:intro-ChatGPT}, and I would like to thank Pavel Etingof and David Penneys for comments and suggestions.

\section{A criterion for rigidity in pivotal braided r-categories}\label{sec:rig-crit}

Let $\Bbbk$ be a field and let $\cC$ be a $\Bbbk$-linear abelian braided monoidal category, with tensor functor $\otimes$, unit object $\vac$, left and right unit isomorphisms $l$ and $r$, associativity isomorphism $\cA$, and braiding isomorphism $c$; see for example \cite{EGNO} for general information on braided monoidal categories. Following \cite{BD}, $K\in\mathrm{Ob}(\cC)$ is \textit{dualizing} if:
\begin{enumerate}

\item For any $X\in\mathrm{Ob}(\cC)$, the functor $\mathrm{Hom}(-\otimes X, K)$ is representable by some $DX\in\mathrm{Ob}(\cC)$, that is, there is a natural isomorphism
\begin{equation}\label{eqn:GV-nat-iso}
\mathrm{Hom}(Y\otimes X, K) \xrightarrow{\sim} \mathrm{Hom}(Y,DX)
\end{equation}
for all $Y\in\mathrm{Ob}(\cC)$.

\item The contravariant functor $X\mapsto DX$ is an anti-equivalence; we call it the \textit{duality functor} of $\cC$.

\end{enumerate}
Then $\cC$ is an \textit{r-category} if $\vac$ is dualizing. From now on, we assume that $\cC$ is an r-category.

For $X\in\mathrm{Ob}(\cC)$, take $Y=DX$ and $K=\vac$ in \eqref{eqn:GV-nat-iso}, and let $e_X: DX\otimes X\rightarrow\vac$ be the morphism mapped to $\Id_{DX}$ under \eqref{eqn:GV-nat-iso}. Then \eqref{eqn:GV-nat-iso} implies that the pair $(DX,e_X)$ is characterized by a universal property: for any $Y\in\mathrm{Ob}(\cC)$ and morphism $f: Y\otimes X\rightarrow\vac$, there is a unique morphism $g: Y\rightarrow DX$ such that the diagram
\begin{equation*}
\xymatrixcolsep{3pc}
 \xymatrix{
 Y\tens X \ar[rd]^{f} \ar[d]_{g\tens\Id_X} & \\
DX\tens X \ar[r]_(.6){e_X} & \vac \\
 }
\end{equation*}
commutes. In particular, there is a unique isomorphism $\varphi:\vac\rightarrow D\vac$ such that the diagram
 \begin{equation}\label{diag:phi_char}
   \xymatrixcolsep{3pc}
   \begin{matrix}
   \xymatrix{
   \vac\otimes \vac \ar[rd]^{l_V=r_V} \ar[d]_{\varphi\tens\Id_\vac} & \\
   D\vac\tens \vac \ar[r]_(.6){e_\vac} & \vac\\}
   \end{matrix}
  \end{equation}
commutes, with inverse $D\vac\xrightarrow{r_\vac^{-1}} D\vac\otimes\vac\xrightarrow{e_\vac}\vac$ (see \cite[Remark 1.1(iv)]{BD}), and for a morphism $f: X\rightarrow Y$, the dual $Df: DY\rightarrow DX$ is the unique morphism such that the diagram
\begin{equation}\label{diag:f'_def}
 \xymatrixcolsep{4pc}
 \begin{matrix}
 \xymatrix{
 DY\otimes X \ar[r]^{\Id_{DY}\otimes f} \ar[d]_{Df\otimes\Id_X} & DY\otimes Y \ar[d]^{e_{Y}} \\
 DX\otimes X \ar[r]_(.57){e_{X}} & \vac \\
 }
 \end{matrix}
\end{equation}
commutes. For brevity and clarity, we will use graphical calculus in this section and the next. For example, we write the commutative diagrams \eqref{diag:phi_char} and \eqref{diag:f'_def} in graphical calculus as:
\begin{equation*}
\begin{matrix}
    \begin{tikzpicture}[scale = 1, baseline = {(current bounding box.center)}, line width=0.75pt] 
    \draw[dashed] (0,0) -- (0,1.25) .. controls (0,2.15) and (1.5,2.15) .. (1.5,1.25) -- (1.5,0);
    \draw[dashed] (.75,1.75) -- (.75,2.5);
    \node at (0,-.25) {$\vac$};
    \node at (1.5,-.25) {$\vac$};
    \node at (.75,2.75) {$\vac$};
    \node at (0,.75) [draw,minimum width=20pt,minimum height=10pt,thick, fill=white] {$\varphi$};
    \node at (.75,1.85) [draw,minimum width=20pt,minimum height=10pt,thick, fill=white] {$e_\vac$};
    \end{tikzpicture}
   \end{matrix} =
   \begin{matrix}
    \begin{tikzpicture}[scale = 1, baseline = {(current bounding box.center)}, line width=0.75pt] 
    \draw[dashed] (0,0) -- (.75,1) -- (1.5,0);
    \draw[dashed] (.75,1) -- (.75,2);
    \node at (0,-.25) {$\vac$};
    \node at (1.5,-.25) {$\vac$};
    \node at (.75,2.25) {$\vac$};
    \end{tikzpicture}
   \end{matrix}\qquad\text{and}\qquad
 \begin{matrix}
    \begin{tikzpicture}[scale = 1, baseline = {(current bounding box.center)}, line width=0.75pt] 
    \draw[white, double=black, line width = 3pt ] (0,0) -- (0,1.25) .. controls (0,2.15) and (1.5,2.15) .. (1.5,1.25) -- (1.5,0);
    \draw[dashed] (.75,1.75) -- (.75,2.5);
    \node at (0,-.25) {$DY$};
    \node at (1.5,-.25) {$X$};
    \node at (.75,2.75) {$\vac$};
    \node at (0,.75) [draw,minimum width=20pt,minimum height=10pt,thick, fill=white] {$Df$};
    \node at (.75,1.85) [draw,minimum width=20pt,minimum height=10pt,thick, fill=white] {$e_{X}$};
    \end{tikzpicture}
   \end{matrix} =
  \begin{matrix}
    \begin{tikzpicture}[scale = 1, baseline = {(current bounding box.center)}, line width=0.75pt] 
    \draw[white, double=black, line width = 3pt ] (0,0) -- (0,1.25) .. controls (0,2.15) and (1.5,2.15) .. (1.5,1.25) -- (1.5,0);
    \draw[dashed] (.75,1.75) -- (.75,2.5);
    \node at (0,-.25) {$DY$};
    \node at (1.5,-.25) {$X$};
    \node at (.75,2.75) {$\vac$};
    \node at (1.5,.75) [draw,minimum width=20pt,minimum height=10pt,thick, fill=white] {$f$};
    \node at (.75,1.85) [draw,minimum width=20pt,minimum height=10pt,thick, fill=white] {$e_{Y}$};
    \end{tikzpicture}
   \end{matrix}.
\end{equation*}

Now \cite[Definition 5.1]{BD} defines a pivotal structure on $\cC$ to be a natural isomorphism
\begin{equation*}
\psi_{X,Y}: \mathrm{Hom}(X\otimes Y,\vac)\xrightarrow{\sim}\mathrm{Hom}(Y\otimes X, \vac)
\end{equation*}
for $X,Y\in\mathrm{Ob}(\cC)$ satisfying certain properties, and by \cite[Proposition 5.7]{BD}, this is equivalent to a natural isomorphism $\delta: \Id_\cC\rightarrow D^2$ which also satisfies suitable properties.
 But since $\cC$ is braided, it is more convenient to characterize pivotal structures in terms of twists, which are natural isomorphisms $\theta:\Id_\cC\rightarrow\Id_\cC$ satisfying the balancing equation
\begin{equation}\label{eqn:balancing}
\theta_{X\otimes Y}= c_{Y,X}\circ c_{X,Y}\circ(\theta_X\otimes\theta_Y)
\end{equation}
for $X,Y\in\mathrm{Ob}(\cC)$; a twist automatically satisfies $\theta_\vac=\Id_\vac$ (see \cite[Remark 6.6(i) and (ii)]{BD}). In particular, by \cite[Proposition 7.1]{BD}, there is a bijection between pivotal structures $\delta$ on $\cC$ and twists $\theta$, where $\delta$ and $\theta$ are related as follows: for any $X\in\mathrm{Ob}(\cC)$, $\delta_X: X\rightarrow D^2X$ is the unique morphism such that the diagram
\begin{equation}\label{diag:delta_def}
  \xymatrixcolsep{4pc}
  \begin{matrix}
  \xymatrix{
  X\otimes DX \ar[d]_{\delta_X\boxtimes\Id_{DX}} \ar[rd]^{\quad\qquad\til{e}_X:=e_X\circ c_{X,DX}\circ(\theta_X\otimes\Id_{DX})} &  \\
D^2X \otimes DX \ar[r]_(.6){e_{DX}} & \vac \\
  }
  \end{matrix}
 \end{equation}
 commutes. That is,
 \begin{equation*}
 \begin{matrix}
    \begin{tikzpicture}[scale = 1, baseline = {(current bounding box.center)}, line width=0.75pt] 
    \draw[white, double=black, line width = 3pt ] (0,0) -- (0,1.25) .. controls (0,2.15) and (1.5,2.15) .. (1.5,1.25) -- (1.5,0);
    \draw[dashed] (.75,1.75) -- (.75,2.5);
    \node at (0,-.25) {$X$};
    \node at (1.5,-.25) {$DX$};
    \node at (.75,2.75) {$\vac$};
    \node at (0,.75) [draw,minimum width=20pt,minimum height=10pt,thick, fill=white] {$\delta_X$};
    \node at (.75,1.75) [draw,minimum width=20pt,minimum height=10pt,thick, fill=white] {$e_{DX}$};
    \end{tikzpicture}
   \end{matrix} =
   \begin{matrix}
    \begin{tikzpicture}[scale = 1, baseline = {(current bounding box.center)}, line width=0.75pt] 
    \draw[white, double=black, line width = 3pt ] (0,0) -- (0,.5) .. controls (0,1.4) and (1.5,1.4) .. (1.5,.5) -- (1.5,0);
    \draw[dashed] (.75,1) -- (.75,2);
    \node at (0,-.25) {$X$};
    \node at (1.5,-.25) {$DX$};
    \node at (.75,2.25) {$\vac$};
    \node at (.75,1.1) [draw,minimum width=20pt,minimum height=10pt,thick, fill=white] {$\til{e}_{X}$};
    \end{tikzpicture}
   \end{matrix}
   =
  \begin{matrix}
    \begin{tikzpicture}[scale = 1, baseline = {(current bounding box.center)}, line width=0.75pt] 
    \draw (1.5,0) -- (1.5,1.15) .. controls (1.5,1.65) and (0,1.4) .. (0,1.9) -- (0,2);
    \draw[white, double=black, line width = 3pt ] (0,0) -- (0,1.15) .. controls (0,1.65) and (1.5,1.4) .. (1.5,1.9) -- (1.5,2);
    \draw (0,2) .. controls (0,2.85) and (1.5,2.85) .. (1.5,2);
    \draw[dashed] (.75,2.55) -- (.75,3.25);
    \node at (0,-.25) {$X$};
    \node at (1.5,-.25) {$DX$};
    \node at (.75,3.5) {$\vac$};
    \node at (0,.65) [draw,minimum width=20pt,minimum height=10pt,thick, fill=white] {$\theta_X$};
    \node at (.75,2.55) [draw,minimum width=20pt,minimum height=10pt,thick, fill=white] {$e_{X}$};
    \end{tikzpicture}
   \end{matrix}
\end{equation*}
in graphical calculus. See Remark 1.1(ii) and Definition 6.7 of \cite{BD} for the definitions of the morphisms used in \cite[Proposition 7.1]{BD}, which are used to derive \eqref{diag:delta_def}. Both morphisms $e_X: DX\otimes X\rightarrow\vac$ and $\til{e}_X: X\otimes DX\rightarrow\vac$ will play an important role in the following.

Now we discuss criteria for $X\in\mathrm{Ob}(\cC)$ to be rigid. Since $\cC$ is braided, it is enough to consider whether $X$ has a left dual, and we would like $DX$ to be a left dual of $X$ with evaluation morphism $e_X$. Thus we need a coevaluation $i_X: \vac\rightarrow X\otimes DX$ such that the first ``zigzag'' composition
\begin{align*}
X\xrightarrow{l_X^{-1}} \vac\otimes X\xrightarrow{i_X\otimes\Id_X} (X\otimes DX)\otimes X\xrightarrow{\cA_{X,DX,X}^{-1}} X\otimes(DX\otimes X)\xrightarrow{\Id_X\otimes e_X} X\otimes\vac\xrightarrow{r_X} X
\end{align*}
equals $\Id_X$. It will then follow automatically that the second zigzag composition
\begin{align*}
DX\xrightarrow{r_{DX}^{-1}} &\, DX\otimes\vac  \xrightarrow{\Id_{DX}\otimes i_X} DX\otimes (X\otimes DX)\\
&\xrightarrow{\cA_{DX,X,DX}} (DX\otimes X)\otimes DX\xrightarrow{e_X\otimes\Id_{DX}} \vac\otimes DX\xrightarrow{l_{DX}} DX
\end{align*}
equals $\Id_{DX}$ (see for example \cite[Lemma 4.1]{CMY3} or \cite[Lemma 2.1]{EP}), so that $(DX, e_X, i_X)$ will be a left dual of $X$.

To find a coevaluation, we use the unique morphism $\Phi_{X,Y}: DY\otimes DX\rightarrow D(X\otimes Y)$, for $X,Y\in\mathrm{Ob}(\cC)$, characterized by the graphical equality
\begin{equation*}
 \begin{matrix}
    \begin{tikzpicture}[scale = 1, baseline = {(current bounding box.center)}, line width=0.75pt]
    \draw (0,0) -- (0,.5) .. controls (0,.8) .. (.5,.8) .. controls (1,.8) .. (1,.5) -- (1,0);
    \draw[dashed] (1.5, 2.25) -- (1.5,3);
    \draw (.5,1) -- (.5,1.6) .. controls (.5,2.6) and (2.5,2.6) .. (2.5,1.6) -- (2.5, 1);
    \draw (2.5,0) -- (2.5,1);
    \node at (0,-.25) {$DY$};
    \node at (1,-.25) {$DX$};
    \node at (2.5,-.25) {$X\otimes Y$};
    \node at (1.5,3.25) {$\vac$};
    \node at (.5,.9) [draw,minimum width=20pt,minimum height=10pt,thick, fill=white] {$\Phi_{X,Y}$};
    \node at (1.5,2.25) [draw,minimum width=20pt,minimum height=10pt,thick, fill=white] {$e_{X\otimes Y}$};
    \end{tikzpicture}
   \end{matrix} =
   \begin{matrix}
    \begin{tikzpicture}[scale = 1, baseline = {(current bounding box.center)}, line width=0.75pt]
     \draw[white, double=black, line width = 3pt ] (1,0) -- (1,.5) .. controls (1,1.4) and (2,1.4) .. (2,0.5) -- (2,0);                      
    \draw[white, double=black, line width = 3pt ] (0,0) -- (0,1.5) .. controls (0,2.5) and (3,2.5) .. (3,1.5) -- (3,0);
    \draw[dashed] (1.5,1) .. controls (1.5,1.4) .. (0,1.5);
    \draw[dashed] (1.5,2.25) -- (1.5,3);
    \node at (0,-.25) {$DY$};
    \node at (1,-.25) {$DX$};
    \node at (2,-.25) {$X$};
    \node at (3,-.25) {$Y$};
    \node at (1.5,3.25) {$\vac$};
    \node at (1.5,1) [draw,minimum height=10pt,thick, fill=white] {$e_X$};
    \node at (1.5,2.25) [draw,minimum width=20pt,minimum height=10pt,thick, fill=white] {$e_Y$};
    \end{tikzpicture}
   \end{matrix}.
\end{equation*}
Taking $Y=DX$, we set
\begin{equation*}
\Phi_X :=\Phi_{X,DX}\circ(\delta_X\otimes\Id_{DX}): X\otimes DX\rightarrow D(X\otimes DX),
\end{equation*}
which is characterized by the graphical equality
\begin{equation}\label{diag:Phi-def}
 \begin{matrix}
    \begin{tikzpicture}[scale = 1, baseline = {(current bounding box.center)}, line width=0.75pt]
    \draw (0,0) -- (0,.5) .. controls (0,.8) .. (.5,.8) .. controls (1,.8) .. (1,.5) -- (1,0);
    \draw[dashed] (1.5, 2.25) -- (1.5,3);
    \draw (.5,1) -- (.5,1.6) .. controls (.5,2.6) and (2.5,2.6) .. (2.5,1.6) -- (2.5, 1);
    \draw (2.5,0) -- (2.5,1);
    \node at (0,-.25) {$X$};
    \node at (1,-.25) {$DX$};
    \node at (2.5,-.25) {$X\otimes DX$};
    \node at (1.5,3.25) {$\vac$};
    \node at (.5,.9) [draw,minimum width=30pt,minimum height=10pt,thick, fill=white] {$\Phi_{X}$};
    \node at (1.5,2.25) [draw,minimum width=20pt,minimum height=10pt,thick, fill=white] {$e_{X\otimes DX}$};
    \end{tikzpicture}
   \end{matrix} =
   \begin{matrix}
    \begin{tikzpicture}[scale = 1, baseline = {(current bounding box.center)}, line width=0.75pt]
     \draw[white, double=black, line width = 3pt ] (1,0) -- (1,.5) .. controls (1,1.4) and (2,1.4) .. (2,0.5) -- (2,0);                      
    \draw[white, double=black, line width = 3pt ] (0,0) -- (0,1.5) .. controls (0,2.5) and (3,2.5) .. (3,1.5) -- (3,0);
    \draw[dashed] (1.5,1) .. controls (1.5,1.4) .. (0,1.5);
    \draw[dashed] (1.5,2.25) -- (1.5,3);
    \node at (0,-.25) {$X$};
    \node at (1,-.25) {$DX$};
    \node at (2,-.25) {$X$};
    \node at (3,-.25) {$DX$};
    \node at (1.5,3.25) {$\vac$};
    \node at (1.5,1) [draw,minimum height=10pt,thick, fill=white] {$e_X$};
    \node at (1.5,2.25) [draw,minimum width=20pt,minimum height=10pt,thick, fill=white] {$\til{e}_X$};
    \end{tikzpicture}
   \end{matrix}.
\end{equation}
Then by \cite[Corollary 3.6]{BD} (see also \cite[Proposition 4.4.5]{CMY2} and \cite[Proposition 3.2]{McR}), $X$ is rigid if and only if $\Phi_X$ is an isomorphism, in which case the coevaluation $i_X: \vac\rightarrow X\otimes DX$ is the composition
\begin{equation}\label{eqn:i_X-def}
\vac \xrightarrow{\varphi} D\vac \xrightarrow{D\til{e}_X} D(X\otimes DX)\xrightarrow{\Phi_X^{-1}} X\otimes DX,
\end{equation}
using the morphisms defined in \eqref{diag:phi_char}, \eqref{diag:f'_def}, and \eqref{diag:delta_def}.

To show that $\Phi_X$ is an isomorphism, we need a criterion from \cite[Section 3]{McR}. Note that both $D\til{e}_X$ and $\Phi_X$ map into $D(X\otimes DX)$, and that a minimum requirement for the coevaluation candidate \eqref{eqn:i_X-def} to be definable is that $\mathrm{Im}\,D\til{e}_X\subseteq\mathrm{Im}\,\Phi_X$. In fact, if this containment holds and $\Phi_X$ is also injective (but not necessarily surjective), then $i_X$ can be defined.
\begin{thm}[\cite{McR}, Theorem 3.4]\label{thm:injective-criterion}
Let $\cC$ be a $\Bbbk$-linear abelian pivotal braided r-category, and let $X\in\mathrm{Ob}(\cC)$. If $\im D\til{e}_X\subseteq\im\Phi_X$, then $\Phi_X$ is injective. 
\end{thm}
\begin{proof}
In \cite[Theorem 3.4]{McR}, it was assumed that $\cC$ is a module category for a self-contragredient vertex operator algebra, but the proof that $\Phi_X$ is injective if $\im D\til{e}_X\subseteq\im\Phi_X$ works for general $\Bbbk$-linear pivotal r-categories with only a few small changes and remarks.

First, the proof uses \cite[Lemma 3.3]{McR}, which involves a certain morphism
\begin{equation*}
\Psi_X: DX\otimes D(X\otimes DX)\longrightarrow DX
\end{equation*}
defined using the inverse twist $\theta_{DX}^{-1}$. Since here we assume only that $\cC$ is a pivotal braided r-category rather than a ribbon r-category, that is, we do not assume $\theta_{DX}= D\theta_X$, we need to replace $\theta_{DX}^{-1}$ with $D\theta_X^{-1}$ in the definition of $\Psi_X$. Specifically, we now define $\Psi_X$ by the graphical equality
\begin{equation*}
  \begin{matrix}
     \begin{tikzpicture}[scale = 1, baseline = {(current bounding box.center)}, line width=0.75pt]
     \draw (1.5,1) .. controls (1.5,1.7) .. (.75,1.8);
     \draw (0,1) .. controls (0,1.7) .. (.75,1.8);
     \draw (.75,2.2) .. controls (.75,3.3) and (3,3.3) .. (3,2.2) -- (3,1);
     \draw[dashed] (1.875,3) -- (1.875,3.8);
    \node at (0,.75) {$DX$};
    \node at (3,.75) {$X$};
    \node at (1.875,4) {$\vac$};
    \node at (1.5,.75) {$D(X\tens DX)$};
    \node at (.75,2) [draw,minimum width=20pt,minimum height=10pt,thick, fill=white] {$\Psi_X$};
    \node at (1.875,3) [draw,minimum width=20pt,minimum height=10pt,thick, fill=white] {$e_X$};   
     \end{tikzpicture}
   \end{matrix}   = 
   \begin{matrix}
    \begin{tikzpicture}[scale = 1, baseline = {(current bounding box.center)}, line width=0.75pt]
    \draw (0,0) -- (0,1);
     \draw (0,1) -- (0,1.25) .. controls (0, 1.75) and (1.5, 1.5) .. (1.5, 2) .. controls (1.5, 2.5) and (3, 2.25) .. (3,2.75);
     \draw[white, double=black, line width = 3pt ] (1.5,0) -- (1.5,1.25) .. controls (1.5, 1.75) and (0,1.5) .. (0,2) -- (0,3.5) .. controls (0,4.6) and (2.25,4.6) .. (2.25,3.5); 
    \draw[white, double=black, line width = 3pt ] (3,0) -- (3,2) .. controls (3,2.5) and (1.5,2.25) .. (1.5,2.75);
    \draw (1.5,2.75) -- (1.5, 3) -- (2.25, 3.5);
    \draw (2.25,3.5) -- (3,3) -- (3,2.75);
    \draw[dashed] (1.125, 4.2) -- (1.125,5);
    \node at (0,-.25) {$DX$};
    \node at (3,-.25) {$X$};
    \node at (1.125,5.2) {$\vac$};
    \node at (1.5,-.25) {$D(X\tens DX)$};
    \node at (0,.8) [draw,minimum width=20pt,minimum height=10pt,thick, fill=white] {$D\theta_{X}^{-1}$};
    \node at (1.125,4.2) [draw,minimum width=20pt,minimum height=10pt,thick, fill=white] {$e_{X\otimes DX}$};
    \end{tikzpicture}
   \end{matrix} .
\end{equation*}
Then the proof of \cite[Lemma 3.3(1)]{McR} goes through, proving
\begin{equation*}
\Psi_X\circ(\Id_{DX}\otimes\Phi_X) = l_{DX}\circ(e_X\otimes\Id_{DX})\circ\cA_{DX,X,DX}
\end{equation*}
as morphisms $DX\otimes(X\otimes DX)\rightarrow DX$. We also remark that the proof of \cite[Lemma 3.3(2)]{McR} assumes that the pair $(X,\til{e}_X)$ satisfies the same universal property as $(D^2X, e_{DX})$. This assumption remains valid in our setting due to the pivotal structure $\delta_X: X\xrightarrow{\sim} D^2X$.

Now using \cite[Lemma 3.3]{McR}, the proof of \cite[Theorem 3.4]{McR} goes through in our setting with two comments. First, \cite[Theorem 3.4]{McR} concludes that $\Phi_X$ is an isomorphism, but in fact the proof shows that $\Phi_X$ is injective and then uses properties of modules for vertex operator algebras to conclude that $\Phi_X$ is also surjective. So here, we only conclude that $\Phi_X$ is injective. Second, the proof of \cite[Theorem 3.4]{McR} uses the assumption that the functor $-\otimes\mathrm{Ker}\,\Phi_X$ is right exact, or more specifically, that $-\otimes\mathrm{Ker}\,\Phi_X$ preserves surjections. This assumption is still valid here because for $X,Y,Z\in\mathrm{Ob}(\cC)$, there is a natural isomorphism
\begin{equation*}
\mathrm{Hom}(X\otimes Y,Z)\xrightarrow{\sim}\mathrm{Hom}(X,D(Y\otimes D^{-1}Z)),
\end{equation*}
that is, $D(Y\otimes D^{-1}Z)$ is an internal hom of $Y$ and $Z$ in $\cC$ (see \cite[Remark 1.1(iii)]{BD}).
\end{proof}

To improve the result of Theorem \ref{thm:injective-criterion}, we need another result from \cite{McR} whose proof goes through unchanged in our setting. In particular, the only properties of the twist $\theta$ used in the proof are its naturality, the balancing equation \eqref{eqn:balancing}, and $\theta_\vac=\Id_\vac$.
\begin{prop}[\cite{McR}, Proposition 3.5]\label{prop:Phi-dual}
In the setting of Theorem \ref{thm:injective-criterion}, $\Phi_X = D\Phi_X\circ\delta_{X\otimes DX}$ for all $X\in\mathrm{Ob}(\cC)$.
\end{prop}

Using Theorem \ref{thm:injective-criterion} and Proposition \ref{prop:Phi-dual} we can obtain our main rigidity criterion:
\begin{thm}\label{thm:rigidity-criterion}
Let $\cC$ be a $\Bbbk$-linear abelian pivotal braided r-category, and let $X\in\mathrm{Ob}(\cC)$. If $\ker\Phi_X\subseteq\ker\til{e}_X$, then $\Phi_X$ is an isomorphism. In particular, $X$ is rigid with (left) dual $DX$, evaluation $e_X$, and coevaluation $i_X=\Phi_X^{-1}\circ D\til{e}_X\circ\varphi$.
\end{thm}

\begin{proof}
Let $k: \ker\Phi_X\rightarrow X\otimes DX$ and $\til{k}: \ker\til{e}_X\rightarrow X\otimes DX$ be the injective kernel morphisms. We are assuming that there is an injective morphism $q: \ker\Phi_X\rightarrow\ker\til{e}_X$ such that $\til{k}\circ q=k$. To apply Theorem \ref{thm:injective-criterion}, we need to show that $\im D\til{e}_X\subseteq\im\Phi_X$.

Since $\cC$ is abelian, we can take $\im\til{e}_X=\coker \til{k}$ and $\im\Phi_X =\coker k$, so we get a diagram
\begin{equation*}
\xymatrixrowsep{1.5pc}
\xymatrix{
& \ker\til{e}_X \ar[d]^{\til{k}} & &\\
\ker\Phi_X \ar[ru]^q \ar[r]_k & X\otimes DX \ar[d]_{\til{c}} \ar[r]^c & \im\Phi_X \ar[r] \ar@{-->}[ld]^{\exists\,!\,r} & 0\\
& \im\til{e}_X \ar[d] &&\\
& 0 &&\\
}
\end{equation*}
with right exact row, right exact column, and commutative triangles, where the existence and uniqueness of $r$ follows from the universal property of the cokernel $(\mathrm{Im}\,\Phi_X,c)$.
%\begin{equation*}
%\til{c}\circ k = \til{c}\circ\til{k}\circ q =0,
%\end{equation*}
%the universal property of cokernels yields a unique morphism $r: \im\Phi_X\rightarrow\im\til{e}_X$ such that $r\circ c=\til{c}$. 
Note $r$ is surjective because the cokernel morphism $\til{c}$ is. 

Now writing $\Phi_X=j\circ c$ and $\til{e}_X=\til{j}\circ\til{c}$, where $j:\im\Phi_X\rightarrow D(X\otimes DX)$ and $\til{j}: \im\til{e}_X\rightarrow D(X\otimes DX)$ are injective, we get a commutative diagram
\begin{equation*}
\xymatrix{
& X\otimes DX \ar[ld]_{\Phi_X} \ar[d]^c \ar[r]^= & X\otimes DX \ar[d]^{\til{c}} \ar[rd]^{\til{e}_X} & \\
D(X\otimes DX) & \im\Phi_X \ar[l]^(.4){j} \ar[r]_r & \im\til{e}_X \ar[r]_(.55){\til{j}} & \vac\\
}
\end{equation*}
Since $D$ is an anti-equivalence and thus exact, dualizing this diagram leads to 
\begin{equation*}
\xymatrix{
& D(X\otimes DX)    & D(X\otimes DX)  \ar[l]_= & \\
D^2(X\otimes DX) \ar[ru]^{D\Phi_X} \ar[r]_(.55){Dj} & D(\im\Phi_X) \ar[u]_{Dc}   & D(\im\til{e}_X) \ar[u]_{D\til{c}} \ar[l]^{Dr}  & D\vac\ar[l]^(.35){D\til{j}} \ar[lu]_{D\til{e}_X}\\
}
\end{equation*}
where $Dc$, $D\til{c}$, and $Dr$ are injective, and $Dj$ and $D\til{j}$ are surjective. 

Using the pivotal structure on $\cC$ and Proposition \ref{prop:Phi-dual}, we may replace $D^2(X\otimes DX)$ in this diagram with $X\otimes DX$, $D\Phi_X$ with $\Phi_X=D\Phi_X\circ\delta_{X\otimes DX}$, and $Dj$ with $Dj\circ\delta_{X\otimes DX}$ (which is still surjective). Thus by the universal property of $(\im\Phi_X, j, c)$, we get a unique morphism $f: \im\Phi_X\rightarrow D(\im\Phi_X)$ such that the diagram
\begin{equation*}
\xymatrixrowsep{1.5pc}
\xymatrix{
X\otimes DX \ar[rd]^c \ar[rdd]_{Dj\circ\delta_{X\otimes DX}} \ar[rr]^{\Phi_X} && D(X\otimes DX)\\
& \im\Phi_X \ar[ru]^j \ar[d]^f & \\
& D(\im\Phi_X) \ar[ruu]_{Dc} & \\
}
\end{equation*}
commutes. Then $f$ is an isomorphism because $j$ is injective and $Dj\circ\delta_{X\otimes DX}$ is surjective. Similarly, writing $D\til{e}_X=i\circ p$ where $p: D\vac\rightarrow\im D\til{e}_X$ is surjective and $i:\im D\til{e}_X\rightarrow D(X\otimes DX)$ is injective, we get a unique isomorphism $\til{f}: \im D\til{e}_X\rightarrow D(\im\til{e}_X)$ such that
\begin{equation*}
D\til{c}\circ\til{f} =i.
\end{equation*}
We can now define the injection 
\begin{equation*}
g:=f^{-1}\circ Dr\circ\til{f}: \im D\til{e}_X\longrightarrow\im\Phi_X,
\end{equation*}
and then
\begin{equation*}
j\circ g = j\circ f^{-1}\circ Dr\circ\til{f} =Dc\circ Dr\circ\til{f} = D\til{c}\circ \til{f} =i.
\end{equation*}
This shows that $\im D\til{e}_X\subseteq\im\Phi_X$.

It now follows from Theorem \ref{thm:injective-criterion} that $\Phi_X$ is injective. Since $D$ is an anti-equivalence and thus exact, $D\Phi_X$ is surjective, and then by Proposition \ref{prop:Phi-dual}, $\Phi_X$ is also surjective. Thus $\Phi_X$ is an isomorphism, and then $X$ is rigid with dual $DX$ by \cite[Corollary 3.6]{BD}.
\end{proof}

To check the rigidity criterion of Theorem \ref{thm:rigidity-criterion} using moderate growth assumptions, it is helpful to formulate the result in contrapositive form:
\begin{cor}\label{cor:non-rigidity-consequence}
In the setting of Theorem \ref{thm:rigidity-criterion}, let $k:\ker\Phi_X\rightarrow X\otimes DX$ be the injective kernel morphism of $\Phi_X$ for some $X\in\mathrm{Ob}(\cC)$. If $X$ is not rigid, then $\til{e}_X\circ k\neq 0$.
\end{cor}
\begin{proof}
If $X$ is not rigid, then by Theorem \ref{thm:rigidity-criterion}, there does not exist an injection $q:\ker \Phi_X\rightarrow\ker\til{e}_X$ such that $k=\til{k}\circ q$, where $\til{k}: \ker\til{e}_X\rightarrow X\otimes DX$ is the injective kernel morphism of $\til{e}_X$. Thus by the universal property of $(\ker\til{e}_X,\til{k})$, the composition $\til{e}_X\circ k$ cannot be $0$.
\end{proof}

\section{Rigidity from moderate growth conditions}\label{sec:rigidity-proof}

We continue to assume $\cC$ is a $\Bbbk$-linear abelian pivotal braided r-category. We will use Corollary \ref{cor:non-rigidity-consequence} together with methods from \cite{EP} to show that rigidity for $X\in\mathrm{Ob}(\cC)$ follows from bounds on the dimensions of certain morphism spaces related to $X$. To simplify notation, we will generally suppress associativity isomorphisms throughout the exposition.

Fix $X\in\mathrm{Ob}(\cC)$ and set $K=\ker\Phi_X$ with injective kernel morphism $k: K\rightarrow X\otimes DX$. We first define linear maps
\begin{align*}
F_n: \mathrm{End}(X^{\otimes n})&\longrightarrow\hom((DX)^{\otimes n}\otimes X^{\otimes n},\vac)\\
F'_n: \mathrm{End}((DX)^{\otimes n})&\longrightarrow\hom((DX)^{\otimes n}\otimes X^{\otimes n},\vac)
\end{align*}
for any $n\in\ZZ_{\geq 1}$ graphically by
\begin{equation*}
F_n(f) =
 \begin{matrix}
    \begin{tikzpicture}[scale = 1, baseline = {(current bounding box.center)}, line width=0.75pt]
        \draw (0,3.25) -- (0,3.75);
        \draw (2,3.25) -- (2,3.75);
        \draw (4,3.25) -- (4,3.75);
        \draw (6,3.25) -- (6,3.75);
        \draw (0,3.75) -- (0,4.25);
        \draw (2,3.75) -- (2,4.25);
        \draw (4,3.75) -- (4,4.25);
        \draw (6,3.75) -- (6,4.25);
        \draw (2,4.25) .. controls (2,5.25) and (4,5.25) .. (4,4.25);
        \draw (0,4.25) .. controls (0,7) and (6,7) .. (6,4.25);
        \draw[dashed] (3,6.25) -- (3,7);
\node at (0,3) {$DX$};
\node at (1,3) {$\cdots$};
\node at (2,3) {$DX$};
\node at (4,3) {$X$};
\node at (5,3) {$\cdots$};
\node at (6,3) {$X$};
\node at (3,5.75) {$\vdots$};
\node at (3,7.25) {$\vac$};
      \node at (5,3.75) [draw,minimum width=65pt,minimum height=10pt,thick, fill=white] {$f$};
      \node at (3,5) [draw,minimum width=20pt,minimum height=10pt,thick, fill=white] {$e_X$};
      \node at (3,6.25) [draw,minimum width=20pt,minimum height=10pt,thick, fill=white] {$e_X$};
        \end{tikzpicture}
   \end{matrix}
   \end{equation*}
   and
   \begin{equation*}
F'_n(g) =
 \begin{matrix}
    \begin{tikzpicture}[scale = 1, baseline = {(current bounding box.center)}, line width=0.75pt]
        \draw (0,3.25) -- (0,3.75);
        \draw (2,3.25) -- (2,3.75);
        \draw (4,3.25) -- (4,3.75);
        \draw (6,3.25) -- (6,3.75);
        \draw (0,3.75) -- (0,4.25);
        \draw (2,3.75) -- (2,4.25);
        \draw (4,3.75) -- (4,4.25);
        \draw (6,3.75) -- (6,4.25);
        \draw (2,4.25) .. controls (2,5.25) and (4,5.25) .. (4,4.25);
        \draw (0,4.25) .. controls (0,7) and (6,7) .. (6,4.25);
        \draw[dashed] (3,6.25) -- (3,7);
\node at (0,3) {$DX$};
\node at (1,3) {$\cdots$};
\node at (2,3) {$DX$};
\node at (4,3) {$X$};
\node at (5,3) {$\cdots$};
\node at (6,3) {$X$};
\node at (3,5.75) {$\vdots$};
\node at (3,7.25) {$\vac$};
      \node at (1,3.75) [draw,minimum width=65pt,minimum height=10pt,thick, fill=white] {$g$};
      \node at (3,5) [draw,minimum width=20pt,minimum height=10pt,thick, fill=white] {$e_X$};
      \node at (3,6.25) [draw,minimum width=20pt,minimum height=10pt,thick, fill=white] {$e_X$};
        \end{tikzpicture}
   \end{matrix}.
\end{equation*}
Note that $F_n$ and $F'_n$ are linear because we assume $\cC$ is a $\Bbbk$-linear monoidal category.
Then similarly as in \cite[Section 2.2]{EP}, we define a bilinear map
\begin{align*}
\phi_n: \mathrm{End}(X^{\otimes n})\times\mathrm{End}(X^{\otimes n}) & \longrightarrow\hom(K^{\otimes n},\vac)
\end{align*}
such that $\phi_n(f,g)$ is the composition
\begin{equation*}
K^{\otimes n}\xrightarrow{k^{\otimes n}} (X\otimes DX)^{\otimes n}\xrightarrow{\mathrm{Braid}(n)} (DX)^{\otimes n}\otimes X^{\otimes n}\xrightarrow{F_n(f\circ g)} \vac,
\end{equation*}
where 
$\mathrm{Braid}(n)$ is defined recursively by $\mathrm{Braid}(0)=\Id_\vac$ and for $n>0$,
\begin{equation*}
\mathrm{Braid}(n)= \begin{matrix}
    \begin{tikzpicture}[scale = 1, baseline = {(current bounding box.center)}, line width=0.75pt]
    \node at (0,-.25) {$X$}; 
    \node at (2,-.25) {$DX$};
    \node at (5,-.25) {$(X\otimes DX)^{\otimes (n-1)}$};
    \node at (0,2.25) {$DX$}; 
    \node at (2,2.25) {$DX^{\otimes (n-1)}$}; 
    \node at (4,2.25) {$X^{\otimes (n-1)}$}; 
    \node at (6,2.25) {$X$}; 
    \draw (2,0) -- (2,1) .. controls (2,1.7) and (0,1.3) .. (0,2);
    \draw (5,0) -- (5,.5);
    \draw (4,.5) -- (4,1) .. controls (4,1.7) and (2,1.3) .. (2,2);
    \draw (6,.5) -- (6,1) .. controls (6,1.7) and (4,1.3) .. (4,2);
    \draw[white, double=black, line width = 3pt ] (0,0) -- (0,1) .. controls (0,1.7) and (6,1.3) .. (6,2);
    \node at (0,.5) [draw,minimum width=20pt,minimum height=10pt,thick, fill=white] {$\theta_{X}$};
    \node at (5,.5) [draw,minimum width=20pt,minimum height=10pt,thick, fill=white] {$\mathrm{Braid}(n-1)$};
    \end{tikzpicture}
    \end{matrix}.
\end{equation*}
Using \eqref{diag:delta_def} and the recursive definition of $\mathrm{Braid}(n)$, the following holds by induction on $n$:
\begin{lem}\label{lem:B(Id,Id)}
For $n\in\ZZ_{\geq 1}$, $\phi_n(\Id_X^{\otimes n},\Id_X^{\otimes n}) =(\til{e}_X\circ k)^{\otimes n}$
up to unit isomorphisms $\vac^{\otimes n}\xrightarrow{\sim}\vac$.
\end{lem}

We are mainly concerned with endomorphisms of $X^{\otimes n}$ and $(DX)^{\otimes n}$ that are constructed from the braiding $c$ on $\cC$. Thus for $1\leq i\leq n-1$, we set
\begin{align*}
\sigma_i &=\Id_X^{\otimes (i-1)}\otimes c_{X,X}\otimes\Id_X^{\otimes (n-i-1)}\in\mathrm{End}(X^{\otimes n}),\\\sigma_i' & =\Id_{DX}^{\otimes (n-i-1)}\otimes c_{DX,DX}\otimes\Id_{DX}^{\otimes (i-1)}\in\mathrm{End}((DX)^{\otimes n}).
\end{align*}
Our first important observation is the following:
\begin{lem}\label{lem:sigma-sigma'}
For $n\in\ZZ_{\geq 1}$ and $1\leq i\leq n-1$, $F_n(\sigma_i) =F'_n(\sigma_{i}')$ and $F_n(\sigma_i^{-1}) =F'_n((\sigma_{i}')^{-1})$.
\end{lem}
\begin{proof}
Note that $\sigma_i^{\pm 1}$ only affects the $i$ and $i+1$ tensor factors of $X^{\otimes n}$, and in the definition of $F_n(\sigma_i^{\pm 1})$, these tensor factors are paired via $e_X$ with the $n-i+1$ and $n-i$ tensor factors of $(DX)^{\otimes n}$. Similarly, $(\sigma_{i}')^{\pm 1}$ only affects the $n-i$ and $n-i+1$ tensor factors of $(DX)^{\otimes n}$, and in the definition of $F'_n((\sigma_{i}')^{\pm 1})$, these tensor factors are paired via $e_X$ with the $i+1$ and $i$ tensor factors of $X^{\otimes n}$. Thus by naturality of unit isomorphisms, it is enough to prove
\begin{equation*}
F_{2}(\sigma_1^{\pm 1}) = F_{2}'((\sigma_1')^{\pm 1}),
\end{equation*}
which has the graphical proof
\begin{align*}
F_{2}(\sigma_1^{\pm 1}) & =  \begin{matrix}
    \begin{tikzpicture}[scale = 1, baseline = {(current bounding box.center)}, line width=0.75pt]
        \draw (1,0) -- (1,1) .. controls (1,1.6) and (2,1.6) .. (2,1) .. controls (2,.3) and (3,.7) .. (3,0);
        \draw[white, double=black, line width = 3pt ] (0,0) -- (0,1) .. controls (0,2.5) and (3,2.5) .. (3,1) .. controls (3,.3) and (2,.7) .. (2,0);
        \draw[dashed] (1.5,1.3) -- (1.5,2.75);
              \node at (1.5,1.3) [draw,minimum width=10pt,minimum height=10pt,thick, fill=white] {$e_X$};
      \node at (1.5,2.1) [draw,minimum width=20pt,minimum height=10pt,thick, fill=white] {$e_X$};
      \node at (1.5,3) {$\vac$};
      \node at (1.5,3) {$\vac$};
      \node at (0,-.25) {$DX$};
      \node at (1,-.25) {$DX$};
      \node at (2,-.25) {$X$};
      \node at (3,-.25) {$X$};
        \end{tikzpicture}
   \end{matrix} =
   \begin{matrix}
    \begin{tikzpicture}[scale = 1, baseline = {(current bounding box.center)}, line width=0.75pt]
        \draw (1,0) -- (1,1) .. controls (1,1.7) and (2,1.3) .. (2,2) .. controls (2,2.6) and (3,2.6) .. (3,2) .. controls (3,1.3) and (2,1.7) .. (2,1) .. controls (2,.3) and (3,.7) .. (3,0);
        \draw[white, double=black, line width = 3pt ] (0,0) -- (0,3) .. controls (0,3.5) and (1,3.5) .. (1,3) -- (1,2) .. controls (1,1.3) and (3,1.7) .. (3,1) .. controls (3,.3) and (2,.7) .. (2,0);
        \draw[dashed] (.5,3.3) -- (.5,3.9);
        \draw[dashed] (2.5,2.3) .. controls (2.5,3) .. (1,3);
              \node at (2.5,2.3) [draw,minimum width=10pt,minimum height=10pt,thick, fill=white] {$e_X$};
      \node at (.5,3.3) [draw,minimum width=10pt,minimum height=10pt,thick, fill=white] {$e_X$};
      \node at (.5,4.15) {$\vac$};
      \node at (0,-.25) {$DX$};
      \node at (1,-.25) {$DX$};
      \node at (2,-.25) {$X$};
      \node at (3,-.25) {$X$};
        \end{tikzpicture}
   \end{matrix} =
   \begin{matrix}
    \begin{tikzpicture}[scale = 1, baseline = {(current bounding box.center)}, line width=0.75pt]
        \draw (1,0) .. controls (1,.7) and (2,.3) .. (2,1) -- (2,2) .. controls (2,2.5) and (3,2.5) .. (3,2) -- (3,0);
        \draw[white, double=black, line width = 3pt ] (0,0) -- (0,1) .. controls (0,1.5) and (1,1.5) .. (1,1) .. controls (1,.3) and (2,.7) .. (2,0);
        \draw[dashed] (2.5,2.3) -- (2.5,2.9);
        \draw[dashed] (.5,1.3) .. controls (.5,2) .. (2,2);
              \node at (2.5,2.3) [draw,minimum width=10pt,minimum height=10pt,thick, fill=white] {$e_X$};
      \node at (.5,1.3) [draw,minimum width=10pt,minimum height=10pt,thick, fill=white] {$e_X$};
      \node at (2.5,3.15) {$\vac$};
      \node at (0,-.25) {$DX$};
      \node at (1,-.25) {$DX$};
      \node at (2,-.25) {$X$};
      \node at (3,-.25) {$X$};
        \end{tikzpicture}
   \end{matrix}\\
   & =
   \begin{matrix}
    \begin{tikzpicture}[scale = 1, baseline = {(current bounding box.center)}, line width=0.75pt]
        \draw (1,0) .. controls (1,.7) and (2,.3) .. (2,1) .. controls (2,1.7) and (0,1.3) .. (0,2) .. controls (0,3.5) and (3,3.5) .. (3,2) -- (3,0);
        \draw[white, double=black, line width = 3pt ] (0,0) -- (0,1) .. controls (0,1.7) and (1,1.3) .. (1,2) .. controls (1,2.5) and (2,2.5) .. (2,2) .. controls (2,1.3) and (1,1.7) .. (1,1) .. controls (1,.3) and (2,.7) .. (2,0);
        \draw[dashed] (1.5,2.3) -- (1.5,3.75);
              \node at (1.5,2.3) [draw,minimum width=10pt,minimum height=10pt,thick, fill=white] {$e_X$};
      \node at (1.5,3.1) [draw,minimum width=10pt,minimum height=10pt,thick, fill=white] {$e_X$};
      \node at (1.5,4) {$\vac$};
      \node at (0,-.25) {$DX$};
      \node at (1,-.25) {$DX$};
      \node at (2,-.25) {$X$};
      \node at (3,-.25) {$X$};
        \end{tikzpicture}
   \end{matrix} = 
    \begin{matrix}
    \begin{tikzpicture}[scale = 1, baseline = {(current bounding box.center)}, line width=0.75pt]
        \draw (1,0) .. controls (1,.7) and (0,.3) .. (0,1) .. controls (0,2.5) and (3,2.5) .. (3,1) -- (3,0);
        \draw[white, double=black, line width = 3pt ] (0,0) .. controls (0,.7) and (1,.3) .. (1,1) .. controls (1,1.5) and (2,1.5) .. (2,1) -- (2,0);
        \draw[dashed] (1.5,1.3) -- (1.5,2.75);
              \node at (1.5,1.3) [draw,minimum width=10pt,minimum height=10pt,thick, fill=white] {$e_X$};
      \node at (1.5,2.1) [draw,minimum width=10pt,minimum height=10pt,thick, fill=white] {$e_X$};
      \node at (1.5,3) {$\vac$};
      \node at (0,-.25) {$DX$};
      \node at (1,-.25) {$DX$};
      \node at (2,-.25) {$X$};
      \node at (3,-.25) {$X$};
        \end{tikzpicture}
   \end{matrix} = F'_{2}((\sigma'_1)^{\pm 1}),
\end{align*}
where over-crossings of strands here can represent either the braiding $c$ or its inverse $c^{-1}$.
\end{proof}

\begin{cor}\label{cor:f'-existence}
For any $f$ in the subgroup $\langle \sigma_1,\ldots,\sigma_{n-1}\rangle\subseteq\mathrm{Aut}(X^{\otimes n})$ generated by the $\sigma_i$, there exists $f'\in\mathrm{Aut}((DX)^{\otimes n})$ such that $F_{n}(f) =F'_{n}(f')$.
\end{cor}
\begin{proof}
Write $f$ as $\sigma_{i_1}^{\varepsilon_1}\cdots\sigma_{i_j}^{\varepsilon_j}$ for $i_1,\ldots,i_j\in\lbrace 1,\ldots, n-1\rbrace$ and $\varepsilon_1,\ldots,\varepsilon_j\in\lbrace\pm 1\rbrace$. We prove by induction on $j$ that we may take $f'=(\sigma_{i_j}')^{\varepsilon_j}\cdots(\sigma_{i_1}')^{\varepsilon_1}$. There is nothing to prove in the base case $j=0$, and for $j>0$, the definitions of $F_{n}$ and $F'_{n}$, Lemma \ref{lem:sigma-sigma'}, and the induction hypothesis show that
\begin{align*}
F_{n}(f) & = F_{n}(\sigma_{i_1}^{\varepsilon_1})\circ(\Id_{DX}^{\otimes n}\otimes\sigma_{i_2}^{\varepsilon_2}\cdots\sigma_{i_j}^{\varepsilon_j}) = F_{n}'((\sigma_{i_1}')^{\varepsilon_1})\circ(\Id_{DX}^{\otimes n}\otimes\sigma_{i_2}^{\varepsilon_2}\cdots\sigma_{i_j}^{\varepsilon_j}) \\
& = F_{n}(\sigma_{i_2}^{\varepsilon_2}\cdots\sigma_{i_j}^{\varepsilon_j})\circ((\sigma_{i_1}')^{\varepsilon_1}\otimes\Id_X^{\otimes n}) =F'_{n}((\sigma_{i_j}')^{\varepsilon_j}\cdots(\sigma_{i_2}')^{\varepsilon_2})\circ((\sigma_{i_1}')^{\varepsilon_1}\otimes\Id_X^{\otimes n})\\
& = F'_{n}(f'),
\end{align*}
as required.
\end{proof}

Now as in \cite[Section 2.1]{EP}, we map the symmetric group $S_n$ into $\mathrm{Aut}(X^{\otimes n})$ as a subset (not as a subgroup if the braiding on $\cC$ is not symmetric). First, $S_n$ embeds as a subset of the braid group $B_n$ on $n$ strands, with generators $b_1,\ldots,b_{n-1}$, as follows: to $\sigma\in S_n$, written as a reduced word in the standard generators $s_i=(i,i+1)$ for $i=1,\ldots, n-1$, we associate the element of $g_\sigma\in B_n$ which is the same word written in the generators $b_i$. More concretely, as in Equation (2) of \cite{EP}, for any $\sigma\in S_n$,
\begin{equation}\label{eqn:f-sigma-decomp}
g_\sigma =w_{n-1} w_{n-2}\cdots w_2 w_1
\end{equation}
for unique
\begin{equation*}
w_j\in\lbrace \Id\rbrace\cup\lbrace b_i\cdots b_{j-1} b_j\mid 1\leq i\leq j\rbrace
\end{equation*}
for $1\leq j\leq n-1$ (in terms of braid diagrams, $b_i\cdots b_{j-1} b_j$ for $i\leq j$ is the braid on $n$ strands that moves the $j+1$ strand underneath strands $i$ through $j$). Specifically, $w_j$ for $1\leq j\leq n-1$ is determined by the insertion sort algorithm, that is, $w_j=b_i\cdots b_{j-1}b_j$ for the unique $i\in\lbrace 1,\ldots, j+1\rbrace$ such that $\sigma(j+1)$ is the $i$th largest element of $\lbrace \sigma(1),\ldots,\sigma(j+1)\rbrace$. Then different elements of $S_n$ yield different tuples $(w_1,\ldots, w_{n-1})$ of elements of $B_n$.

For $m\leq n$ we embed $B_m\subseteq B_n$ as the subgroup generated by $b_1,\ldots, b_{m-1}$. We recall the following key lemma from \cite{EP}:
\begin{lem}[\cite{EP}, Lemma 2.3]\label{lem:EP}
If $\sigma,\tau\in S_n$ are distinct, then there exist $m\in\lbrace 1,\ldots, n-1\rbrace$, $u,v\in B_m$, and $\varepsilon\in\lbrace \pm 1\rbrace$ such that
\begin{equation*}
g_\sigma^{-1} g_\tau = u b_m^\varepsilon v.
\end{equation*}
\end{lem}

\begin{proof}
As in \eqref{eqn:f-sigma-decomp}, we write $g_\sigma= w_{n-1} w_{n-2}\cdots w_1$ and $g_\tau=v_{n-1} v_{n-2}\cdots v_1$ where
\begin{equation*}
v_j,w_j\in\lbrace\Id\rbrace\cup\lbrace b_i\cdots b_{j-1}b_j\mid 1\leq i\leq j\rbrace
\end{equation*}
for $1\leq j\leq n-1$. Since $\sigma\neq\tau$, let $m\in\lbrace 1,\ldots, n-1\rbrace$ be maximal such that $w_m\neq v_m$. Then
\begin{equation*}
g_\sigma^{-1} g_\tau = w_1^{-1}\cdots w_{m-1}^{-1}(w_m^{-1}v_m)v_{m-1}\cdots v_1,
\end{equation*}
where
\begin{equation*}
w_m^{-1} v_m = b_m^{-1}b_{m-1}^{-1}\cdots b_i^{-1}b_j\cdots b_{m-1}b_m
\end{equation*}
for $i\neq j$. The quick braid diagram calculation in the proof of \cite[Lemma 2.3]{EP} shows that $w_m^{-1} v_m =u' b_m^\varepsilon v'$ for some $u',v'\in B_m$ and $\varepsilon\in\lbrace\pm 1\rbrace$, so then 
\begin{equation*}
g_\sigma^{-1}g_\tau =u b_m^\varepsilon v
\end{equation*}
where $u= w_1^{-1}\cdots w_{m-1}^{-1} u'$ and $v= v'v_{m-1}\cdots v_1$.
\end{proof}

For $\sigma\in S_n$, let $f_\sigma\in\mathrm{Aut}(X^{\otimes n})$ be the image of $g_\sigma$ under
the group homomorphism $B_n\rightarrow\mathrm{Aut}(X^{\otimes n})$ which sends $b_i$ to $\sigma_i$ for $1\leq i\leq n-1$. We can now prove the main theorem:
\begin{thm}\label{thm:main-theorem}
Let $\cC$ be a $\Bbbk$-linear abelian pivotal braided r-category with a simple unit object $\vac$, and let $X\in\mathrm{Ob}(\cC)$. If there exists $n\in\ZZ_{\geq 1}$ such that either
\begin{equation*}
\dim\mathrm{Hom}((DX)^{\otimes n}\otimes X^{\otimes n},\vac)<n!
\end{equation*}
or
\begin{equation*}
\dim\mathrm{End}(X^{\otimes n})<n!,
\end{equation*}
then $X$ is rigid with left dual $DX$.
\end{thm}
\begin{proof}
We will show that if $X$ is not rigid, then the morphisms $F_{n}(f_\sigma)\in\mathrm{Hom}((DX)^{\otimes n}\otimes X^{\otimes n},\vac)$ and the morphisms $f_\sigma\in\mathrm{End}(X^{\otimes n})$ are linearly independent for any $n\in\ZZ_{\geq 1}$, and thus both morphism spaces have dimension at least $n!$.

Thus assume $X$ is not rigid; by Corollary \ref{cor:non-rigidity-consequence}, $\til{e}_X\circ k\neq 0$, where $k:\mathrm{Ker}\,\Phi_X\rightarrow X\otimes DX$ is the injective kernel morphism of $\Phi_X: X\otimes DX\rightarrow D(X\otimes DX)$ defined in \eqref{diag:Phi-def}, and $\til{e}_X$ is defined in \eqref{diag:delta_def}. Since $\vac$ is simple, this implies $\til{e}_X\circ k$ is surjective. Moreover, as observed in the proof of Theorem \ref{thm:injective-criterion}, internal homs exist in $\cC$ by \cite[Remark 1.1(iii)]{BD}, and this implies that the tensor product $\otimes$ is right exact and thus preserves surjections. Thus $(\til{e}_X\circ k)^{\otimes n}$ is surjective for any $n\in\ZZ_{\geq 1}$.

Now similar to the proof of \cite[Theorem 1.1]{EP}, the key claim is that the bilinear form $\phi_{n}$ defined in the paragraph before Lemma \ref{lem:B(Id,Id)} satisfies
\begin{equation}\label{eqn:claim}
\phi_{n}(f_\sigma^{-1},f_\tau) =\delta_{\sigma,\tau}(\til{e}_X\circ k)^{\otimes n}
\end{equation}
for all $\sigma,\tau\in S_n$, up to unit isomorphisms $\vac^{\otimes n}\xrightarrow{\sim} \vac$. Before proving the claim, we show how it implies that the morphisms $f_\sigma\in\mathrm{End}(X^{\otimes n})$ for $\sigma\in S_n$ are linearly independent. If there is a relation $\sum_{\tau\in S_n} a_\tau f_\tau=0$ for $a_\tau\in\Bbbk$, then bilinearity of $\phi_{n}$ and the claim imply that
\begin{equation*}
0 =\phi_{n}\bigg(f_\sigma^{-1}, \sum_{\tau\in S_n} a_\tau f_\tau\bigg) =\sum_{\tau\in S_n} a_\tau \phi_{n}(f_\sigma^{-1},f_\tau)=a_\sigma(\til{e}_X\circ k)^{\otimes n}
\end{equation*}
for any $\sigma\in S_n$. Since $(\til{e}_X\circ k)^{\otimes n}$ is surjective and thus non-zero, it follows that $a_\sigma=0$ for all $\sigma\in S_n$, proving that the endomorphisms $f_\sigma$ are linearly independent.

Similarly, the claim implies that the morphisms $F_{n}(f_\sigma)\in\mathrm{Hom}((DX)^{\otimes n}\otimes X^{\otimes n},\vac)$ for $\sigma\in S_n$ are linearly independent. Indeed, if there is a relation $\sum_{\tau\in S_n} a_\tau F_{n}(f_\tau)=0$ for $a_\tau\in\Bbbk$, then using Corollary \ref{cor:f'-existence} and the definitions of $F_{n}$ and $F'_{n}$, we have
\begin{align*}
0 & =\bigg(\sum_{\tau\in S_n} a_\tau F_{n}(f_\tau)\bigg)\circ((f_\sigma^{-1})'\otimes\Id_X^{\otimes n}) =\sum_{\tau\in S_n} a_\tau F_{n}(f_\tau)\circ((f_\sigma^{-1})'\otimes\Id_X^{\otimes n})\\
& =\sum_{\tau\in S_n} a_\tau F'_{n}((f_\sigma^{-1})')\circ(\Id_{DX}^{\otimes n}\otimes f_\tau) = \sum_{\tau\in S_n} a_\tau F_{n}(f_\sigma^{-1})\circ(\Id_{DX}^{\otimes n}\otimes f_\tau)\\
& =\sum_{\tau\in S_n} a_\tau F_{n}(f_\sigma^{-1}\circ f_\tau).
\end{align*}
 for any $\sigma\in S_n$. Composing both sides with $\mathrm{Braid}(n)\circ k^{\otimes n}$ yields
\begin{align*}
0=\bigg(\sum_{\tau\in S_n} a_\tau F_n(f_\sigma^{-1}\circ f_\tau)\bigg)\circ\mathrm{Braid}(n)\circ k^{\otimes n} = \sum_{\tau\in S_n} a_\tau \phi_{n}(f_\sigma^{-1},f_\tau) = a_\sigma (\til{e}_X\circ k)^{\otimes n}
\end{align*}
for all $\sigma\in S_n$, so again $a_\sigma=0$ for all $\sigma\in S_n$ because $(\til{e}_X\circ k)^{\otimes n}\neq 0$.

It remains to prove \eqref{eqn:claim}. The case $\sigma=\tau$ is covered by Lemma \ref{lem:B(Id,Id)}, so assume $\sigma\neq\tau$. By Lemma \ref{lem:EP}, $f_\sigma^{-1}\circ f_\tau=u\circ\sigma_m^\varepsilon\circ v$ for some $m\in\lbrace 1,\ldots, n-1\rbrace$, $u,v\in\langle \sigma_1,\ldots,\sigma_{m-1}\rangle$, and $\varepsilon\in\lbrace\pm 1\rbrace$. Then by the definition of $\phi_{n}$, the recursive definition of $\mathrm{Braid}(n)$, naturality of the braiding isomorphisms, and \eqref{diag:delta_def}, 
\begin{equation*}
\phi_{n}(f_\sigma^{-1},f_\tau) = (\til{e}_X\circ k)^{\otimes (n-m-1)}\otimes\left[F_{m+1}(u\circ \sigma_m^\varepsilon\circ v)\circ\mathrm{Braid}(m+1)\circ k^{\otimes (m+1)} \right].
\end{equation*}
So it is enough to show that the expression in brackets is $0$. In fact, by the definitions,
\begin{align*}
F_{m+1}(u\circ\sigma_m^\varepsilon\circ v)\circ\mathrm{Braid}(m+1) =\begin{matrix}
    \begin{tikzpicture}[scale = 1, baseline = {(current bounding box.center)}, line width=0.75pt]
    \node at (0,-.25) {$X$};
    \node at (1.5,-.25) {$DX$};
    \node at (3,-.25) {$X$}; 
    \node at (4.5,-.25) {$DX$};
    \node at (6.75,-.25) {$(X\otimes DX)^{\otimes (m-1)}$};
    \node at (0,2.25) {$X$};
    \node at (1.5,2.25) {$DX$};
    \node at (3,2.25) {$DX$}; 
    \node at (4.5,2.25) {$(DX)^{\otimes (m-1)}$}; 
    \node at (6.25,2.25) {$X^{\otimes (m-1)}$}; 
    \node at (7.5,2.25) {$X$}; 
    \draw (4.5,0) -- (4.5,1) .. controls (4.5,1.7) and (3,1.3) .. (3,2);
    \draw (6.75,0) -- (6.75,.5);
    \draw (6,.5) -- (6,1) .. controls (6,1.7) and (4.5,1.3) .. (4.5,2);
    \draw (7.5,.5) -- (7.5,1) .. controls (7.5,1.7) and (6,1.3) .. (6,2);
    \draw[white, double=black, line width = 3pt ] (3,0) -- (3,1) .. controls (3,1.7) and (7.5,1.3) .. (7.5,2);
    \draw (0,0) -- (0,2);
    \draw (1.5,0) -- (1.5,2);
    \draw (1.5,2.5) .. controls (1.5,3.2) and (0,2.8) .. (0,3.5) -- (0,6) .. controls (0,7.5) and (7.5,7.5) .. (7.5,6) -- (7.5,3.5);
    \draw (3, 2.5) .. controls (3,3.2) and (1.5,2.8) .. (1.5,3.5) -- (1.5,6) -- (6,6) -- (6,3.5) .. controls (6,2.8) and (7.5,3.2) .. (7.5,2.5);
    \draw (4.5,2.5) .. controls (4.5,3.2) and (3,2.8) .. (3,3.5) -- (3,6) -- (4.5,6) -- (4.5,3.5) .. controls (4.5,2.8) and (6,3.2) .. (6,2.5);
    \draw[dashed] (3.75,6) -- (3.75,7.75);
     \draw[white, double=black, line width = 3pt ] (0,2.5) .. controls (0,3.2) and (7.5,2.8) .. (7.5,3.5);
    \node at (0,1) [draw,minimum width=20pt,minimum height=10pt,thick, fill=white] {$\theta_{X}$};
    \node at (3,.5) [draw,minimum width=20pt,minimum height=10pt,thick, fill=white] {$\theta_{X}$};
    \node at (6.75,.5) [draw,minimum width=20pt,minimum height=10pt,thick, fill=white] {$\mathrm{Braid}(m-1)$};
    \node at (5.25,4) [draw,minimum width=50pt,minimum height=10pt,thick, fill=white] {$v$};
    \node at (6.75,5) [draw,minimum width=50pt,minimum height=10pt,thick, fill=white] {$c_{X,X}^{\varepsilon}$};
     \node at (3.75,6) [draw,minimum width=140pt,minimum height=10pt,thick, fill=white] {$F_{m}(u)$};
      \node at (3.75,7.15) [draw,minimum width=20pt,minimum height=10pt,thick, fill=white] {$e_X$};
     \node at (3.75,8) {$\vac$};
    \end{tikzpicture}
    \end{matrix}.
\end{align*}
By Corollary \ref{cor:f'-existence},
\begin{equation*}
F_{m}(u) =F'_{m}(u') =\begin{matrix}
    \begin{tikzpicture}[scale = 1, baseline = {(current bounding box.center)}, line width=0.75pt]
        \draw (0,0) -- (0,1.5) .. controls (0,3) and (6,3) .. (6,1.5) -- (6,0);
        \draw (2,0) -- (2,1.5) -- (4,1.5) -- (4,0);
        \draw[dashed] (3,1.5) -- (3,3.25);
              \node at (3,2.6) [draw,minimum width=10pt,minimum height=10pt,thick, fill=white] {$e_X$};
      \node at (3,1.5) [draw,minimum width=10pt,minimum height=10pt,thick, fill=white] {$F_{m-1}(\Id_X^{\otimes (m-1)})$};
      \node at (1,.5) [draw,minimum width=65pt,minimum height=10pt,thick, fill=white] {$u'$};
      \node at (3,3.5) {$\vac$};
      \node at (0,-.25) {$DX$};
      \node at (2,-.25) {$(DX)^{\otimes(m-1)}$};
      \node at (4,-.25) {$X^{\otimes(m-1)}$};
      \node at (6,-.25) {$X$};
        \end{tikzpicture}
   \end{matrix}.
\end{equation*}
We put this back into the composition for $F_{m+1}(u\circ\sigma_m^\varepsilon\circ v)\circ\mathrm{Braid}(m+1)$ and use naturality of the braiding isomorphisms to move both $v$ and $(\Id_{DX}\otimes F_{m-1}(\Id_X^{\otimes m}))\circ(u'\otimes\Id_X^{\otimes (m-1)})$ through the braiding of the first factor of $X$ with $DX\otimes(X\otimes DX)^{\otimes m}$. The result is:
\begin{equation}\label{thm:main-thm-calc1}
F_{m+1}(u\circ\sigma_m^\varepsilon\circ v)\circ\mathrm{Braid}(m+1) =\begin{matrix}
    \begin{tikzpicture}[scale = 1, baseline = {(current bounding box.center)}, line width=0.75pt]
    \node at (0,-.25) {$X$};
    \node at (2,-.25) {$DX$};
    \node at (5,-.25) {$(X\otimes DX)^{\otimes m}$};
    \node at (0,1.25) {$X$};
    \node at (2,1.25) {$DX$};
    \node at (4,1.25) {$DX$}; 
    \node at (6,1.25) {$X$};
    \draw (0,0) -- (0,1);
    \draw (2,0) -- (2,1);
    \draw (5,0) -- (5,.5);
    \draw (4,.5) -- (4,1);
    \draw (6,.5) -- (6,1);
    \draw[dashed] (3,4) -- (3,5.6);
    \draw (4,1.5) .. controls (4,2.2) and (2,1.8) .. (2,2.5) -- (2,3.5) .. controls (2,4.2) and (4,4.2) .. (4,3.5) -- (4,2.5) .. controls (4,1.8) and (6,2.2) .. (6,1.5);
    \draw (2,1.5) .. controls (2,2.2) and (0,1.8) .. (0,2.5) -- (0,3.5) .. controls (0,5.5) and (6,5.5) .. (6,3.5) -- (6,2.5);
     \draw[white, double=black, line width = 3pt ] (0,1.5) .. controls (0,2.2) and (6,1.8) .. (6,2.5);
    \node at (0,.5) [draw,minimum width=20pt,minimum height=10pt,thick, fill=white] {$\theta_{X}$};
    \node at (5,.5) [draw,minimum width=60pt,minimum height=10pt,thick, fill=white] {$G$};
    \node at (5,3) [draw,minimum width=60pt,minimum height=10pt,thick, fill=white] {$c_{X,X}^{\varepsilon}$};
     \node at (3,4) [draw,minimum width=20pt,minimum height=10pt,thick, fill=white] {$e_X$};
      \node at (3,5) [draw,minimum width=20pt,minimum height=10pt,thick, fill=white] {$e_X$};
     \node at (3,5.8) {$\vac$};
    \end{tikzpicture}
    \end{matrix}
    \end{equation}
    for some $G: (X\otimes DX)^{\otimes m}\rightarrow DX\otimes X$ which is no longer important. We analyze the remaining braidings graphically using Lemma \ref{lem:sigma-sigma'}, naturality of braiding, \eqref{diag:delta_def}, and \eqref{diag:Phi-def}:
    \begin{align}\label{eqn:main-thm-calc2}
    \begin{matrix}
    \begin{tikzpicture}[scale = 1, baseline = {(current bounding box.center)}, line width=0.75pt]
    \node at (0,.25) {$X$};
    \node at (1,.25) {$DX$};
    \node at (2,.25) {$DX$}; 
    \node at (3,.25) {$X$};
    \draw[dashed] (1.5,4) -- (1.5,5.25);
    \draw (2,.5) -- (2,1.5) .. controls (2,2.2) and (1,1.8) .. (1,2.5) -- (1,3.5) .. controls (1,4.2) and (2,4.2) .. (2,3.5) -- (2,2.5) .. controls (2,1.8) and (3,2.2) .. (3,1.5) -- (3,.5);
    \draw (1,.5) -- (1,1.5) .. controls (1,2.2) and (0,1.8) .. (0,2.5) -- (0,3.5) .. controls (0,5) and (3,5) .. (3,3.5) -- (3,2.5);
     \draw[white, double=black, line width = 3pt ] (0,.5) -- (0,1.5) .. controls (0,2.2) and (3,1.8) .. (3,2.5);
      \node at (0,1) [draw,minimum width=20pt,minimum height=10pt,thick, fill=white] {$\theta_X$};
    \node at (2.5,3) [draw,minimum width=32pt,minimum height=10pt,thick, fill=white] {$c_{X,X}^{\varepsilon}$};
     \node at (1.5,3.9) [draw,minimum width=20pt,minimum height=10pt,thick, fill=white] {$e_X$};
      \node at (1.5,4.6) [draw,minimum width=20pt,minimum height=10pt,thick, fill=white] {$e_X$};
     \node at (1.5,5.4) {$\vac$};
    \end{tikzpicture}
    \end{matrix} & =
    \begin{matrix}
    \begin{tikzpicture}[scale = 1, baseline = {(current bounding box.center)}, line width=0.75pt]
    \node at (0,1.25) {$X$};
    \node at (1,1.25) {$DX$};
    \node at (2,1.25) {$DX$}; 
    \node at (3,1.25) {$X$};
    \draw[dashed] (1.5,4.25) -- (1.5,5.75);
  \draw (2,1.5) -- (2,2.5) .. controls (2,3.2) and (1,2.8) .. (1,3.5) -- (1,3.75) .. controls (1,4.45) and (2,4.45) .. (2,3.75) -- (2,3.5) .. controls (2,2.8) and (3,3.2) .. (3,2.5) -- (3,1.5);
    \draw (1,1.5) -- (1,2.5) .. controls (1,3.2) and (0,2.8) .. (0,3.5) -- (0,4) .. controls (0,5.5) and (3,5.5) .. (3,4) -- (3,3.5);
     \draw[white, double=black, line width = 3pt ] (0,1.5) -- (0,2.5) .. controls (0,3.2) and (3,2.8) .. (3,3.5);
    \node at (1.5,2) [draw,minimum width=32pt,minimum height=10pt,thick, fill=white] {$c_{DX,DX}^{\varepsilon}$};
        \node at (0,2) [draw,minimum width=20pt,minimum height=10pt,thick, fill=white] {$\theta_X$};
     \node at (1.5,4.15) [draw,minimum width=20pt,minimum height=10pt,thick, fill=white] {$e_X$};
      \node at (1.5,5.1) [draw,minimum width=20pt,minimum height=10pt,thick, fill=white] {$e_X$};
     \node at (1.5,5.9) {$\vac$};
    \end{tikzpicture}
    \end{matrix} =
    \begin{matrix}
    \begin{tikzpicture}[scale = 1, baseline = {(current bounding box.center)}, line width=0.75pt]
    \node at (0,1.25) {$X$};
    \node at (1,1.25) {$DX$};
    \node at (2,1.25) {$DX$}; 
    \node at (3,1.25) {$X$};
    \draw[dashed] (.5,3.9) -- (.5,4.65);
    \draw[dashed] (2.5,2.9) .. controls (2.5,3.5) .. (1,3.5);
  \draw (2,1.5) -- (2,2.5) .. controls (2,3.2) and (3,3.2) .. (3,2.5) -- (3,1.5);
    \draw (1,1.5) -- (1,3.5) .. controls (1,4.2) and (0,4.2) .. (0,3.5) -- (0,1.5);
    \node at (1.5,2) [draw,minimum width=32pt,minimum height=10pt,thick, fill=white] {$c_{DX,DX}^{\varepsilon}$};
     \node at (2.5,2.9) [draw,minimum width=20pt,minimum height=10pt,thick, fill=white] {$e_X$};
      \node at (.5,3.9) [draw,minimum width=20pt,minimum height=10pt,thick, fill=white] {$\til{e}_X$};
     \node at (.5,4.85) {$\vac$};
    \end{tikzpicture}
    \end{matrix}\nonumber\\
    & =
    \begin{matrix}
    \begin{tikzpicture}[scale = 1, baseline = {(current bounding box.center)}, line width=0.75pt]
    \node at (0,1.25) {$X$};
    \node at (1,1.25) {$DX$};
    \node at (2,1.25) {$DX$}; 
    \node at (3,1.25) {$X$};
    \draw[dashed] (1.5,2.9) -- (1.5,4.65);
  \draw (0,1.5) -- (0,2.75) .. controls (0,4.25) and (3,4.25) .. (3,2.75) -- (3,1.5);
    \draw (1,1.5) -- (1,2.5) .. controls (1,3.2) and (2,3.2) .. (2,2.5) -- (2,1.5);
    \node at (2.5,2) [draw,minimum width=32pt,minimum height=10pt,thick, fill=white] {$c^{-\varepsilon}$};
     \node at (1.5,2.9) [draw,minimum width=20pt,minimum height=10pt,thick, fill=white] {$e_X$};
      \node at (1.5,3.85) [draw,minimum width=20pt,minimum height=10pt,thick, fill=white] {$\til{e}_X$};
     \node at (1.5,4.8) {$\vac$};
    \end{tikzpicture}
    \end{matrix} =
    \begin{matrix}
    \begin{tikzpicture}[scale = 1, baseline = {(current bounding box.center)}, line width=0.75pt]
    \node at (0,1.25) {$X$};
    \node at (1,1.25) {$DX$};
    \node at (2,1.25) {$DX$}; 
    \node at (3,1.25) {$X$};
    \draw (0,1.5) -- (0,2);
    \draw (1,1.5) -- (1,2);
    \draw (2,1.5) -- (2,2);
    \draw (3,1.5) -- (3,2);
    \draw[dashed] (1.5,3) -- (1.5,3.75);
  \draw (.5, 2) -- (.5,2.5) .. controls (.5, 3.5) and (2.5,3.5) .. (2.5,2.5) -- (2.5, 2);
    \node at (2.5,2) [draw,minimum width=32pt,minimum height=10pt,thick, fill=white] {$c^{-\varepsilon}$};
     \node at (.5,2) [draw,minimum width=32pt,minimum height=10pt,thick, fill=white] {$\Phi_X$};
      \node at (1.5,3.1) [draw,minimum width=20pt,minimum height=10pt,thick, fill=white] {$e_{X\otimes DX}$};
     \node at (1.5,4) {$\vac$};
    \end{tikzpicture}
    \end{matrix}.
    \end{align}
    Combining \eqref{thm:main-thm-calc1} and \eqref{eqn:main-thm-calc2}, we get
\begin{equation*}
F_{m+1}(u\circ\sigma^\varepsilon\circ v)\circ\mathrm{Braid}(m+1)\circ k^{\otimes (m+1)} =    \begin{matrix}
    \begin{tikzpicture}[scale = 1, baseline = {(current bounding box.center)}, line width=0.75pt]
    \node at (.5,-.25) {$K$};
    \node at (2.5,-.25) {$K^{\otimes m}$};
    \node at (0,1.25) {$X$};
    \node at (1,1.25) {$DX$};
    \node at (2,1.25) {$DX$}; 
    \node at (3,1.25) {$X$};
    \draw (.5,0) -- (.5,.5);
    \draw (2.5,0) -- (2.5,.5);
    \draw (0,1.5) -- (0,2);
    \draw (1,1.5) -- (1,2);
    \draw (2,1.5) -- (2,2);
    \draw (3,1.5) -- (3,2);
     \draw (0,.5) -- (0,1);
    \draw (1,.5) -- (1,1);
    \draw (2,.5) -- (2,1);
    \draw (3,.5) -- (3,1);
    \draw[dashed] (1.5,3) -- (1.5,3.75);
  \draw (.5, 2) -- (.5,2.5) .. controls (.5, 3.5) and (2.5,3.5) .. (2.5,2.5) -- (2.5, 2);
     \node at (2.5,.5) [draw,minimum width=32pt,minimum height=10pt,thick, fill=white] {$G\circ k^{\otimes m}$};
     \node at (.5,.5) [draw,minimum width=32pt,minimum height=10pt,thick, fill=white] {$k$};
    \node at (2.5,2) [draw,minimum width=32pt,minimum height=10pt,thick, fill=white] {$c^{-\varepsilon}$};
     \node at (.5,2) [draw,minimum width=32pt,minimum height=10pt,thick, fill=white] {$\Phi_X$};
      \node at (1.5,3.1) [draw,minimum width=20pt,minimum height=10pt,thick, fill=white] {$e_{X\otimes DX}$};
     \node at (1.5,4) {$\vac$};
    \end{tikzpicture}
    \end{matrix} =0
\end{equation*}
since $\Phi_X\circ k =0$. This completes the proof of the theorem.
\end{proof}

We present two corollaries of the main theorem:
\begin{cor}\label{cor:loc-fin}
In the setting of Theorem \ref{thm:main-theorem}, suppose $\cC$ is locally finite abelian and for any simple $X\in\mathrm{Ob}(\cC)$, there exists $n\in\ZZ_{\geq 1}$ such that either
\begin{equation*}
\dim\mathrm{Hom}((DX)^{\otimes n}\otimes X^{\otimes n},\vac)<n!
\end{equation*}
or
\begin{equation*}
\dim\mathrm{End}(X^{\otimes n})<n!.
\end{equation*}
Then $\cC$ is rigid, with left duals given by the r-category duality functor $D$.
\end{cor}
\begin{proof}
By Theorem \ref{thm:main-theorem}, every simple $X\in\mathrm{Ob}(\cC)$ is rigid with left dual $DX$. Then because $\cC$ is locally finite, every object of $\cC$ has finite length, and thus $\cC$ is rigid by \cite[Theorem 3.13]{CMSY} (see also \cite[Theorem 4.4.1]{CMY2}), with lefts duals given by $D$.
\end{proof}

\begin{thm}\label{thm:fin-abelian-rigid}
Let $\cC$ be a locally finite $\Bbbk$-linear abelian pivotal braided r-category with simple unit object $\vac$ and finitely many isomorphism classes of simple objects. Then $\cC$ is rigid with left duals given by the r-category duality functor $D$.
\end{thm}
\begin{proof}
If $X\in\mathrm{Ob}(\cC)$ is simple, then so is $DX$ because $D$ is an anti-equivalence. Thus by Corollary \ref{cor:loc-fin}, it suffices show that there exist constants $A,B>0$ such that for $N\in\ZZ_{\geq 1}$,
\begin{equation*}
\dim\mathrm{Hom}(X_0\otimes X_1\otimes\cdots\otimes X_N,\vac)\leq A\cdot B^N
\end{equation*}
for any simple $X_0,X_1,\ldots,X_N\in\mathrm{Ob}(\cC)$. Indeed, let $\lbrace X^{(i)}\rbrace_{i=1}^M$ be a set of representatives of the isomorphism classes of simple objects of $\cC$.
Since $\cC$ is a locally finite abelian category, all morphism spaces are finite dimensional and every $X\in\mathrm{Ob}(\mathcal{C})$ has finite Jordan-H\"{o}lder length $\ell(X)$. Thus we can set
\begin{equation*}
A =\dim\mathrm{End}(\vac),\qquad B=\max_{1\leq i,j\leq M} \ell(X^{(i)}\otimes X^{(j)}).
\end{equation*}
 Then it is sufficient to prove:
\begin{enumerate}

\item $\dim\mathrm{Hom}(X,\vac)\leq A\cdot\ell(X)$ for any $X\in\mathrm{Ob}(\cC)$, and

\item $\ell(X_0\otimes X_1\otimes\cdots\otimes X_N)\leq B^N$ for any $N\in\ZZ_{\geq 0}$.

\end{enumerate}

We prove (1) by induction on $\ell(X)$, with nothing to prove in the base case $X=0$. If $\ell(X)>0$, then there is a short exact sequence
\begin{equation*}
0\longrightarrow Y\longrightarrow X\longrightarrow Z\longrightarrow 0
\end{equation*}
such that $Y$ is simple and $\ell(Z)=\ell(X)-1$. Since Hom functors are left exact, 
\begin{equation*}
0\longrightarrow\mathrm{Hom}(Z,\vac)\longrightarrow\mathrm{Hom}(X,\vac)\longrightarrow\mathrm{Hom}(Y,\vac)
\end{equation*}
is left exact and thus
\begin{align*}
\dim\mathrm{Hom}(X,\vac) & \leq\dim\mathrm{Hom}(Y,\vac)+\dim\mathrm{Hom}(Z,\vac)\\
&\leq A+A\cdot\ell(Z) =A\cdot(\ell(Z)+1)=A\cdot\ell(X),
\end{align*}
as required, using the definition of $A$ and the induction hypothesis.

Similarly, we prove (2) by induction on $N$, with nothing to prove in the base case $N=0$. If $N>0$, then set $L=\ell(X_1\otimes\cdots\otimes X_N)$, so $L\leq B^{N-1}$ by the induction hypothesis. There is a filtration
\begin{equation*}
0 = Y_0\subseteq Y_1\subseteq\cdots\subseteq Y_{L-1}\subseteq Y_L=X_1\otimes\cdots\otimes X_N
\end{equation*} 
by subobjects such that $Y_l/Y_{l-1}$ is simple for $l=1,\ldots, L$. We now prove by induction on $l$ that $\ell(X_0\otimes Y_l)\leq B\cdot l$ for all $l\in\lbrace 0, 1,\ldots, L\rbrace$. The base case $l=0$ is obvious, and if $l>0$, then because $\otimes$ is right exact (again recall \cite[Remark 1.1(iii)]{BD}),
\begin{equation*}
X_0\otimes Y_{l-1}\longrightarrow X_0\otimes Y_l\longrightarrow X_0\otimes(Y_l/Y_{l-1})\longrightarrow 0
\end{equation*}
is right exact. It follows that
\begin{align*}
\ell(X_0\otimes Y_l)\leq \ell(X_0\otimes Y_{l-1})+\ell(X_0\otimes(Y_l/Y_{l-1}))\leq B\cdot (l-1)+B = B\cdot l
\end{align*}
as required, using the induction hypothesis and the definition of $B$. In particular,
\begin{equation*}
\ell(X_0\otimes X_1\otimes\cdots\otimes X_N) =\ell(X_0\otimes Y_L)\leq B\cdot L\leq B^N,
\end{equation*}
as required. This completes the proof of the theorem.
\end{proof}

\section{Results for vertex operator algebras}\label{sec:VOA-app}

One of the main motivations for Theorem \ref{thm:main-theorem} is its applications to vertex operator algebras, which are rich sources of $\CC$-linear braided monoidal categories. In particular, Huang showed recently in \cite{Hu}, based on the logarithmic vertex tensor category theory developed by Huang, Lepowsky, and Zhang in \cite{HLZ1}-\cite{HLZ8}, that the category of $C_1$-cofinite modules for any vertex operator algebra is a braided monoidal category. However, it is usually difficult to show directly that a module for a vertex operator algebra is rigid, and thus the non-constructive proof of Theorem \ref{thm:main-theorem} will be useful for general results.

\subsection{General theorems}

A vertex operator algebra is a graded vector space $V=\bigoplus_{n\in\ZZ} V_{(n)}$ over $\CC$, such that $\dim V_{(n)}<\infty$ for all $n$ and $V_{(n)}=0$ for sufficiently negative $n$, equipped with a vertex operator map
\begin{align*}
Y_V: V & \longrightarrow \mathrm{End}_\CC(V)[[x,x^{-1}]]\\
v & \longmapsto Y_V(v,x) =\sum_{n\in\ZZ} v_n\,x^{-n-1},
\end{align*}
a vacuum vector $\mathbf{1}\in V_{(0)}$, and a conformal vector $\omega\in V_{(2)}$, satisfying several axioms specified, for example, in \cite[Section 3.1]{LL}. In particular, the coefficients of
\begin{equation*}
Y_V(\omega,x)=\sum_{n\in\ZZ} L_n\,x^{-n-2}
\end{equation*}
(that is, $L_n=\omega_{n+1}$) generate a representation of the Virasoro Lie algebra on $V$. Also, $V_{(n)}$ is the $L_0$-eigenspace with eigenvalue $n$; if $v\in V_{(n)}$, we write $n=\mathrm{wt}\,v$ and say that $n$ is the conformal weight of $v$.

The definition of vertex operator algebra can be relaxed to allow $V$ to be $\frac{1}{2}\ZZ$-graded instead of $\ZZ$-graded. We also get the definition of vertex operator superalgebra by adding a $\ZZ/2\ZZ$-grading $V=V^0\oplus V^1$ compatible with the $\ZZ$- or $\frac{1}{2}\ZZ$-grading, requiring all structure maps to be even, and adding appropriate sign factors in the axioms. In this paper, we also always assume that a vertex operator superalgebra is either $\ZZ$-graded or is $\frac{1}{2}\ZZ$-graded such that $V^i=\bigoplus_{n\in\frac{i}{2}+\ZZ} V_{(n)}$ for $i=0,1$. This will guarantee that braided monoidal categories of modules for vertex operator superalgebras have a ribbon twist.

A grading-restricted generalized module
for a vertex operator (super)algebra $V$ (or $V$-module for short) is a $\CC$-graded (super)space $W=\bigoplus_{h\in\CC} W_{[h]}$, such that $\dim W_{[h]}<\infty$ for all $h\in\CC$ and $W_{[h]}=0$ for $\mathrm{Re}\,h$ sufficiently negative, equipped with a vertex operator map
\begin{align*}
Y_W: V & \longrightarrow \mathrm{End}_\CC(W)[[x,x^{-1}]]\\
v & \longmapsto Y_W(v,x) =\sum_{n\in\ZZ} v_n\,x^{-n-1}
\end{align*}
satisfying several axioms specified, for example, in \cite[Definition 3.1]{CKM} (which in particular covers the superalgebra generality). Again, the coefficients $L_n$ of $Y_W(\omega,x)$ generate a representation of the Virasoro algebra on $W$, and for $h\in\CC$, $W_{[h]}$ is the generalized $L_0$-eigenspace with eigenvalue $h$.

If $V$ is a $\ZZ$-graded vertex operator (super)algebra, then any $V$-module $W$ has a natural $\ZZ_{\geq 0}$-grading $W=\bigoplus_{n=0}^\infty W(n)$ such that
\begin{equation*}
v_m W(n)\subseteq W(n+\mathrm{wt}\,v-m-1)
\end{equation*}
for $v\in V$, $m\in\ZZ$, and $n\in\ZZ_{\geq 0}$, defined as follows. Let $\lbrace\eta_i\rbrace_{i\in I}$ by the set of cosets $\eta\in\CC/\ZZ$ such that $W_{[h]}\neq 0$ for some $h\in\eta$, and let $h_i\in\eta_i$ have minimal real part such that $W_{[h_i]}\neq 0$ ($h_i$ exists because $W_{[h]}=0$ for $\mathrm{Re}\,h$ sufficiently negative). Then
\begin{equation*}
W=\bigoplus_{h\in\CC} W_{[h]} =\bigoplus_{i\in I}\bigoplus_{n=0}^\infty W_{[h_i+n]} =\bigoplus_{n=0}^\infty\bigoplus_{i\in I} W_{[h_i+n]},
\end{equation*}
so we can take $W(n)=\bigoplus_{i\in I} W_{[h_i+n]}$. For each $i\in I$, the subspace $\bigoplus_{n=0}^\infty W_{[h_i+n]}$ is a $V$-module direct summand of $W$. So if $W$ is finitely generated as a $V$-module, then $I$ is finite, and thus $\dim W(n)<\infty$ for all $n\in\ZZ_{\geq 0}$. If $V$ is $\frac{1}{2}\ZZ$-graded, then any $V$-module has an analogous $\frac{1}{2}\ZZ_{\geq 0}$-grading.

If $V$ is a $\ZZ$- or $\frac{1}{2}\ZZ$-graded vertex operator (super)algebra and $W$ is a $V$-module, then the graded dual (super)space $W'=\bigoplus_{h\in\CC} W_{[h]}^*$ is also a $V$-module, with $Y_{W'}$ defined by
\begin{equation*}
\langle Y_{W'}(v,x)w', w\rangle =(-1)^{\vert v\vert\vert w'\vert}\langle w', Y_W(e^{xL_1}(-x^{-2})^{L_0} v,x^{-1}w\rangle
\end{equation*}
for $v\in V$, $w\in W$, and $w'\in W'$ \cite{FHL}, where $\vert\,\,\,\vert$ denotes parity in the superalgebra case. If $V$ is $\frac{1}{2}\ZZ$-graded, then this formula defines two contragredient module structures, depending on the choice of square root of $-1$, but these two module structures are isomorphic (see the discussion preceding \cite[Remark 3.5]{CKM}).

A $V$-module $W$ is \textit{$C_1$-cofinite} if $\dim W/C_1(W)<\infty$, where
\begin{equation*}
C_1(W) =\mathrm{span}\lbrace v_{-1} w\mid w\in W, v\in V, \mathrm{wt}\,v>0\rbrace.
\end{equation*}
Any $C_1$-cofinite module $W$ is finitely generated as a $V$-module, for example by basis elements of a graded subspace $X$ such that $W=C_1(W)\oplus X$ as a graded vector space (see for example \cite[Proposition 2.1]{CMY1}).
Thus if $W$ is $C_1$-cofinite, the homogeneous subspaces $W(n)$ in the $\ZZ_{\geq 0}$- or $\frac{1}{2}\ZZ_{\geq 0}$-grading discussed above are finite dimensional.
Let $\cC^1_V$ be the category of $C_1$-cofinite grading-restricted generalized $V$-modules; the following theorem summarizes the tensor-categorical properties of $\cC^1_V$ that are already in the literature:
\begin{thm}\label{thm:old-C1-properties}
Let $V$ be a $\ZZ$-graded vertex operator algebra, or a $\ZZ$- or $\frac{1}{2}\ZZ$-graded vertex operator superalgebra, such that $V\cong V'$ as a $V$-module, and assume that $\cC^1_V$ is closed under contragredient modules. Then $\cC^1_V$, equipped with the braided monoidal structure of \cite{HLZ8}, is a locally finite $\CC$-linear abelian ribbon (in particular, pivotal braided) r-category, with duality functor $D$ given by taking contragredients.
\end{thm}
\begin{proof}
By \cite[Theorem 1.2]{Hu} (and the paragraph in \cite{Hu} following its statement), $\cC^1_V$ is a $\CC$-linear braided monoidal category as in \cite{HLZ8}, with unit object $V$, and $\cC^1_V$ is locally finite abelian by \cite[Theorem 2.9]{McR-Deligne}. Since $V\cong V'$ as a $V$-module, results from \cite[Section 5]{FHL}, or \cite{HLZ2} easily imply that $V$ is a dualizing object of $\cC^1_V$ with duality functor $D$ given by taking contragredients (see for example \cite[Section 5.1]{CKM2}, \cite[Section 4.3]{CMY2}, or \cite[Theorem 2.12]{ALSW}). Thus $\cC^1_V$ is a braided r-category, and it is also ribbon (and thus pivotal) with twist $\theta_W=e^{2\pi i L_0}$ if $V$ is $\ZZ$-graded or $\theta_W=P_We^{2\pi i L_0}$ (where $P_W$ is the parity involution of the vector superspace $W$) if $V$ is a $\frac{1}{2}\ZZ$-graded superalgebra such that $V^i=\bigoplus_{n\in\frac{i}{2}+\ZZ} V_{(n)}$ for $i=0,1$. 
\end{proof}

To apply Theorem \ref{thm:main-theorem} to the category $\cC^1_V$ of $C_1$-cofinite $V$-modules, we need two lemmas in order to check the first moderate growth condition in Theorem \ref{thm:main-theorem}.
\begin{lem}\label{lem:C1-est-hom}
Let $V$ be a simple vertex operator (super)algebra and let $N\in\frac{1}{2}\ZZ_{\geq 0}$ be maximal such that $V_{(-N)}\neq 0$. If $W$ is a $C_1$-cofinite grading-restricted generalized $V$-module, then 
\begin{equation*}
\dim \mathrm{Hom}(W,V)\leq(\dim W/C_1(W))(\dim V_{(-N)}).
\end{equation*}
\end{lem}
\begin{proof}
Let $W=\bigoplus_{n\in\frac{1}{2}\NN} W(n)$ be the $\frac{1}{2}\ZZ_{\geq 0}$-grading of $W$ discussed above (or $\ZZ_{\geq 0}$-grading if $V$ is $\ZZ$-graded). Since $W$ is $C_1$-cofinite, each $W(n)$ is finite dimensional. Also, $C_1(W)$ is a $\frac{1}{2}\ZZ_{\geq 0}$-graded subspace of $W$ because $v_m W(n)\subseteq W(n+\mathrm{wt}\,v-m-1)$ for $v\in V$, $m\in\ZZ$, $n\in\frac{1}{2}\ZZ_{\geq 0}$. As explained in the proof of \cite[Theorem 3.7]{CY}, there is a finite filtration
\begin{equation*}
0= W_0\subseteq W_1\subseteq W_2\subseteq\cdots\subseteq W_{M-1}\subseteq W_M=W
\end{equation*}
such that $W_i/W_{i-1}$ is generated as a $V$-module by $(W_i/W_{i-1})(0)$ for $i\in\lbrace 1,\ldots, M\rbrace$. Specifically, assuming that $W_{i-1}$ has already been chosen, we take $W_i$ to be the submodule of $W$ such that $W_i/W_{i-1}$ is the submodule of $W/W_{i-1}$ generated by $(W/W_{i-1})(0)$. This process terminates after finitely many steps because $W$ is $C_1$-cofinite and thus finitely generated.

For $i\in\lbrace 1,\ldots, M\rbrace$, let $X_i$ be a finite-dimensional subspace of $W$ such that 
\begin{equation*}
X_i\cap W_{i-1}=0,\qquad (X_i+W_{i-1})/W_{i-1} = (W/W_{i-1})(0).
\end{equation*}
 By the construction of $W_i$, we can take $X_i\subseteq W(n_i)$ for $n_i\in\frac{1}{2}\ZZ_{\geq 0}$ such that 
\begin{equation*}
0 = n_1<n_2<\cdots <n_{M-1}<n_M,
\end{equation*}
and $W$ is generated as a $V$-module by $X=\bigoplus_{i=1}^M X_i$. Moreover, $\bigoplus_{n<n_i} W(n)\subseteq W_{i-1}$ for $1\leq i\leq M$. We claim that $X\cap C_1(W)=0$, which will show that $\dim X\leq \dim W/C_1(W)$.

To prove this claim, since both $X$ and $C_1(W)$ are $\frac{1}{2}\ZZ_{\geq 0}$-graded, it is enough to show that $X_i\cap C_1(W)=0$ for each $i\in\lbrace 1,\ldots, M\rbrace$. So suppose $x=\sum_{j=1}^J v^{(j)}_{-1} w^{(j)}\in X_i\cap C_1(W)$ for some $i$ and for $w^{(j)}\in W$ and $v^{(j)}\in V$ such that $\mathrm{wt}\,v^{(j)}>0$. Then we may assume
\begin{equation*}
w^{(j)}\in W(n_i-\mathrm{wt}\,v^{(j)})
\end{equation*}
for all $j$, which implies $w^{(j)}\in\bigoplus_{n<n_i} W(n)\subseteq W_{i-1}$ since $\mathrm{wt}\,v^{(j)}>0$, and thus $x\in W_{i-1}$. This means $x=0$ since $X_i\cap W_{i-1}=0$ by assumption, proving the claim. In particular, $\dim X\leq\dim W/C_1(W)$.

Next, we prove by induction on $i$ that $\dim\mathrm{Hom}(W_i,V)\leq\sum_{j=1}^i \dim\mathrm{Hom}(W_j/W_{j-1}, V)$ for $1\leq i\leq M$. There is nothing to prove for the base case $i=1$ since $W_0=0$, so now suppose $i>1$. Then there is an exact sequence
\begin{equation*}
0\longrightarrow W_{i-1}\longrightarrow W_i\longrightarrow W_i/W_{i-1}\longrightarrow 0,
\end{equation*}
so because Hom functors are left exact, there is a left exact sequence
\begin{equation*}
0\longrightarrow\mathrm{Hom}(W_i/W_{i-1}, V)\longrightarrow\mathrm{Hom}(W_i,V)\longrightarrow\mathrm{Hom}(W_{i-1},V).
\end{equation*}
Thus using the induction hypothesis,
\begin{align*}
\dim\mathrm{Hom}(W_i,V) & \leq \dim\hom(W_{i-1},V)+\dim\hom(W_i/W_{i-1},V)\\
&\leq\sum_{j=1}^i\dim\hom(W_j/W_{j-1},V),
\end{align*}
as required. In particular, $\dim\hom(W,V)\leq\sum_{i=1}^M \dim\mathrm{Hom}(W_i/W_{i-1},V)$.

Finally, we claim that $\dim\hom(W_i/W_{i-1},V)\leq(\dim X_i)(\dim V_{(-N)})$. To prove this, first note that the linear map
\begin{align*}
\hom(W_i/W_{i-1},V) &\longrightarrow \hom_\CC((W_i/W_{i-1})(0),V)\\
f & \longmapsto f\vert_{(W_i/W_{i-1})(0)}
\end{align*}
is injective since $W_i/W_{i-1}$ is generated by $(W_i/W_{i-1})(0)$ as a $V$-module. Second, for any $f\in\hom(W_i/W_{i-1},V)$ and $x\in (W_i/W_{i-1})(0)$, we have
\begin{equation*}
v_n f(x) = f(v_n x) =0
\end{equation*}
for $v\in V$ and $n\in\ZZ$ such that $\mathrm{wt}\,v-n-1<0$. Thus
\begin{equation*}
f((W_i/W_{i-1})(0))\subseteq S:=\lbrace s\in V\mid v_n s=0\,\,\text{if}\,\,\mathrm{wt}\,v-n-1<0\rbrace.
\end{equation*}
We want to show that this $\frac{1}{2}\ZZ_{\geq 0}$-graded subset $S$ equals $V_{(-N)}$. Indeed, $V_{(-N)}\subseteq S$ because $V_{(-N)}=V(0)$ in the $\frac{1}{2}\ZZ_{\geq 0}$-grading of $V$ as a $V$-module. For the reverse inclusion, suppose $v\in V_{(n)}$ for $n>-N$ and $s\in V_{(-N)}$ are both non-zero. Because $V$ is simple, \cite[Proposition 4.5.6]{LL} shows that $s=\sum_j u^{(j)}_{m_j} v$ for certain $u^{(j)}\in V$ and $m_j\in\ZZ$; we may assume $u^{(j)}$ is homogeneous and $u^{(j)}_{m_j} v\neq 0$ for each $j$, which implies that $v\notin S$. That is, $S\cap V_{(n)}=0$ for $n>-N$, and thus $S=V_{(-N)}$ as desired.
 Then $f\mapsto f\vert_{(W_i/W_{i-1})(0)}$ gives a linear injection $\hom(W_i/W_{i-1},V)\rightarrow\hom_\CC((W_i/W_{i-1})(0),V_{(-N)})$, and thus
\begin{align*}
\dim\hom(W_i/W_{i-1},V) & \leq\dim\hom_\CC((W_i/W_{i-1})(0),V_{(-N)})\\
& =\dim\hom_\CC(X_i,V_{(-N)})=(\dim X_i)(\dim V_{(-N)})
\end{align*}
because $X_i\cong (W_i/W_{i-1})(0)$ as vector spaces. This proves the claim.

Now we put all the claims together and conclude
\begin{align*}
\dim\hom(W,V) & \leq\sum_{i=1}^M \dim\hom(W_i/W_{i-1},V)\leq\sum_{i=1}^M (\dim X_i)(\dim V_{(-N)})\\
& =(\dim X)(\dim V_{(-N)})\leq (\dim W/C_1(W))(\dim V_{(-N)}),
\end{align*}
which proves the lemma.
\end{proof}

Following standard vertex operator algebra notation, we use $\boxtimes$ to denote the tensor product operation on the category $\cC^1_V$ of $C_1$-cofinite $V$-modules. The next lemma was first discovered in a conformal field theory setting by Nahm \cite{Na}. In the vertex algebraic context, it follows from the proof of \cite[Key Theorem]{Mi} and is also proved explicitly in \cite[Theorem 3.5]{Hu}:
\begin{lem}\label{lem:C1-estimate-tens}
If $W_1$ and $W_2$ are $C_1$-cofinite grading-restricted generalized modules for any vertex operator (super)algebra, then
\begin{equation*}
\dim W_1\boxtimes W_2/C_1(W_1\boxtimes W_2)\leq (\dim W_1/C_1(W_1))(\dim W_2/C_2(W_2)).
\end{equation*}
\end{lem}

We can now state and prove the main application to vertex operator algebras:
\begin{thm}\label{thm:C1-cofinite-rigidity}
Let $V$ be a simple $\ZZ$-graded vertex operator algebra or superalgebra, or a simple $\frac{1}{2}\ZZ$-graded vertex operator superalgebra such that $V^i=\bigoplus_{n\in\frac{i}{2}+\ZZ} V_{(n)}$ for $i=0,1$. Assume that $V\cong V'$ as a $V$-module and that the category $\cC_V^1$ of $C_1$-cofinite grading-restricted generalized $V$-modules is closed under contragredient modules. Then $\cC_V^1$, equipped with the braided monoidal structure of \cite{HLZ8}, is a locally finite $\CC$-linear abelian braided ribbon category. In particular, $\cC^1_V$ is rigid with duals given by contragredient modules.
\end{thm}
\begin{proof}
By Theorem \ref{thm:old-C1-properties}, $\cC^1_V$ equipped with the braided monoidal structure of \cite{HLZ8} is already a locally finite $\CC$-linear abelian ribbon r-category, with duality functor $D$ given by taking contragredients. It only remains to show that $\cC^1_V$ is rigid. Let $N\in\frac{1}{2}\ZZ_{\geq 0}$ be maximal such that $V_{(-N)}\neq 0$. Then if $W\in\mathrm{Ob}(\cC_V^1)$,
\begin{align*}
\dim\hom((W')^{\boxtimes n}\boxtimes W^{\boxtimes n},V) &\leq(\dim V_{(-N)})(\dim\, (W')^{\boxtimes n}\boxtimes W^{\boxtimes n}/C_1((W')^{\boxtimes n}\boxtimes W^{\boxtimes n}))\\
& \leq (\dim V_{(-N)})(\dim W'/C_1(W'))^n(\dim W/C_1(W))^n
\end{align*}
for any $n\in\ZZ_{\geq 1}$ by Lemmas \ref{lem:C1-est-hom} and \ref{lem:C1-estimate-tens}. So
\begin{equation*}
\dim \hom((W')^{\boxtimes n}\boxtimes W^{\boxtimes n},V)< n!
\end{equation*}
for sufficiently large $n\in\ZZ_{\geq 1}$, and thus $W$ is rigid with dual $W'$ by Theorem \ref{thm:main-theorem}. Since $W$ was arbitrary, $\cC_V^1$ is rigid with duals given by contragredient modules.
\end{proof}

A vertex operator algebra is \textit{$C_2$-cofinite} if $\dim V/C_2(V) <\infty$, where 
\begin{equation*}
C_2(V) =\mathrm{span}\lbrace u_{-2}v\mid u,v\in V\rbrace.
\end{equation*}
It is shown in \cite{Hu-C2} that if $V$ is $C_2$-cofinite and $\ZZ_{\geq 0}$-graded with $V_{(0)}=\CC\vac$, then the category $\rep(V)$ of all grading-restricted generalized $V$-modules is a finite abelian category and also admits the braided monoidal structure of \cite{HLZ8}. The assumption in \cite{Hu-C2} that $V_{(0)}=\CC\vac$ is only needed to guarantee existence of certain spanning sets of $V$ and of finitely-generated $V$-modules, which were constructed by Gaberdiel-Neitzke \cite{GN} and Buhl \cite{Bu} under the assumption that $V_{(0)}=\CC\vac$. But such spanning sets exist for general $\ZZ_{\geq 0}$-graded $C_2$-cofinite vertex operator algebras by \cite[Lemma 2.4]{Mi-C2}, so in fact $\mathrm{Rep}(V)$ is a finite abelian category and a braided monoidal category if $V$ is any $\ZZ_{\geq 0}$-graded $C_2$-cofinite vertex operator algebra. The spanning set result \cite[Lemma 2.4]{Mi-C2} can also be used to show that $\rep(V)=\cC_V^1$ in this case. It was conjectured in \cite{Hu-conj} that for a suitable self-contragredient $C_2$-cofinite vertex operator algebra $V$, $\rep(V)$ should be rigid and further a (not necessarily semisimple) modular tensor category. We can now use Theorem \ref{thm:main-theorem} and its consequences to complete the proof of this conjecture:
\begin{thm}\label{thm:C2-cofinite-rigidity}
If $V$ is a simple $\ZZ_{\geq 0}$-graded self-contragredient $C_2$-cofinite vertex operator algebra, then the category $\rep(V)$ of grading-restricted generalized $V$-modules is a (not necessarily semisimple) modular tensor category, with braided monoidal structure as in \cite{HLZ8} and duals given by contragredient modules.
\end{thm}
\begin{proof}
From the discussion preceding the theorem, and from the assumptions of the theorem, $\rep(V)$ is a finite $\CC$-linear abelian ribbon (in particular, pivotal braided) r-category with simple unit object $V$, with braided monoidal structure as in \cite{HLZ8} and duality functor $D$ given by contragredients. In particular, $\rep(V)$ has finitely many isomorphism classes of simple objects, so $\rep(V)$ is rigid by Theorem \ref{thm:fin-abelian-rigid}. Alternatively, $\rep(V)$ is rigid by Theorem \ref{thm:C1-cofinite-rigidity} because $\rep(V)=\cC_V^1$. Then $\rep(V)$ is a (not necessarily semisimple) modular tensor category by \cite[Main Theorem 1]{McR}.
\end{proof}

\subsection{Examples}

There are many vertex operator algebras which are known to satisfy the conditions of Theorems \ref{thm:C1-cofinite-rigidity} and \ref{thm:C2-cofinite-rigidity}. Although the conclusion of Theorem \ref{thm:C2-cofinite-rigidity} has been proved previously for many specific examples of $C_2$-cofinite vertex operator algebras, many more examples are expected to exist, including general finite group orbifolds of known $C_2$-cofinite vertex operator algebras and Feigin-Tipunin vertex operator algebras \cite{FT, Su}.

 Here we highlight the application of Theorem \ref{thm:C1-cofinite-rigidity} to a few families for which $\cC_V^1$ has not yet been proved to be rigid in all cases:
\begin{exam}
It is shown in \cite{CJORY} that if $V$ is any Virasoro vertex operator algebra at any central charge, then $\cC_V^1$ is closed under contragredient modules. Moreover, in this case, $V$ is self-contragredient if it is simple, so the category of $C_1$-cofinite grading-restricted generalized modules for any simple Virasoro vertex operator algebra is a braided ribbon category by Theorem \ref{thm:C1-cofinite-rigidity}.
\end{exam}

\begin{exam}
It is shown in \cite{CMOY} that if $V$ is any vertex operator superalgebra associated to the $N=1$ super Virasoro algebra at any central charge, then $\cC_V^1$ is closed under contragredient modules. Moreover, in this case, $V$ is self-contragredient if it is simple, so the category of $C_1$-cofinite grading-restricted generalized modules for any simple $N=1$ super Virasoro vertex operator superalgebra is a braided ribbon category by Theorem \ref{thm:C1-cofinite-rigidity}.
\end{exam}

Finally, we consider Kazhdan-Lusztig categories of modules for affine Lie algebras. Let $\mathfrak{g}$ be a finite-dimensional simple Lie algebra over $\CC$ with re-scaled Killing form $\langle\,,\,\rangle$ such that $\langle\alpha,\alpha\rangle=2$ for long roots $\alpha$, and let
\begin{equation*}
\widehat{\mathfrak{g}}=\mathfrak{g}\otimes\CC[t,t^{-1}]\oplus\CC\mathbf{k}
\end{equation*}
be the affine Lie algebra associated to $\mathfrak{g}$, where $\mathbf{k}$ is central and
\begin{equation*}
[a\otimes t^m, b\otimes t^n] =[a,b]\otimes t^{m+n}+m\langle a,b\rangle\delta_{m+n,0}\mathbf{k} 
\end{equation*}
for $a,b\in\mathfrak{g}$, $m,n\in\ZZ$. A $\widehat{\mathfrak{g}}$-module has level $k\in\CC$ if $\mathbf{k}$ acts by the scalar $k$. For any $k\in\CC$ and $\mathfrak{g}$-module $M$, the Weyl module $W^k_M$ of level $k$ is the induced $\widehat{\mathfrak{g}}$-module
\begin{equation*}
W^k_M=U(\widehat{\mathfrak{g}})\otimes_{U(\widehat{\mathfrak{g}}_{\geq 0})} M,
\end{equation*}
where $\widehat{\mathfrak{g}}_{\geq 0}=\mathfrak{g}\otimes\CC[t]\oplus\CC\mathbf{k}$ and for $x\in M$,
\begin{equation*}
(a\otimes t^n)\cdot x =\delta_{n,0} (a\cdot x),\qquad\mathbf{k}\cdot x=kx
\end{equation*}
for $a\in\mathfrak{g}$, $n\in\ZZ_{\geq 0}$. If $k\neq -h^\vee$, where $h^\vee$ is the dual Coxeter number of $\mathfrak{g}$, and if $M=M_\lambda$ is the simple finite-dimensional $\mathfrak{g}$-module of highest weight $\lambda$, then $W^k_{M_\lambda}$ has a unique simple $\widehat{\mathfrak{g}}$-module quotient, which we denote by $L^k_\lambda$.

Following \cite[Definition 2.15]{KL1}, we define the Kazhdan-Lusztig category $KL^k(\mathfrak{g})$ to be the category of finite-length $\widehat{\mathfrak{g}}$-modules of level $k$ whose Jordan-H\"{o}lder composition factors come from the $L^k_\lambda$ for dominant integral weights $\lambda$ of $\mathfrak{g}$. In \cite{KL1}-\cite{KL4}, Kazhdan and Lusztig showed that if $\mathfrak{g}$ is simply laced and $k+h^\vee\notin\QQ_{\geq 0}$ (with some possible exceptions for $k+h^\vee\in\QQ_{<0}$ with small numerator), then $KL^k(\mathfrak{g})$ is a braided ribbon category equivalent to the category of finite-dimensional weight modules for a suitable version of the quantum group of $\mathfrak{g}$ at parameter $q=e^{\pi i/(k+h^\vee)}$. 

We can now use Theorem \ref{thm:C1-cofinite-rigidity} to give an easier proof that $KL^k(\mathfrak{g})$ is rigid, and at the same time we remove the simply laced assumption and the assumption that the numerator of $\vert k+h^\vee\vert$ is sufficiently large if $k+h^\vee\in\QQ_{<0}$:

\begin{thm}
Let $\mathfrak{g}$ be any finite-dimensional simple Lie algebra over $\CC$ and take $k\in\CC$ such that $k+h^\vee\notin\QQ_{\geq 0}$. Then $KL^k(\mathfrak{g})$ is the category of $C_1$-cofinite grading-restricted generalized modules for the universal affine vertex operator algebra of $\mathfrak{g}$ at level $k$ constructed in \cite{FZ}, and $KL^k(\mathfrak{g})$ equipped with the braided monoidal structure of \cite{HLZ8} is a braided ribbon category.
\end{thm}
\begin{proof}
Let $V^k(\mathfrak{g})$ be the universal affine vertex operator algebra of $\mathfrak{g}$ at level $k$ constructed in \cite{FZ}; it is a vertex operator algebra structure on the Weyl module $W^k_{M_0}$ induced from the trivial one-dimensional $\mathfrak{g}$-module $\CC\vac$, such that
\begin{equation*}
Y_{V^k(\mathfrak{g})}((a\otimes t^{-1})\cdot\vac,x) =\sum_{n\in\ZZ} (a\otimes t^n)\,x^{-n-1}
\end{equation*}
for $a\in\mathfrak{g}$.
The proof that $KL^k(\mathfrak{g})=\cC^1_{V^k(\mathfrak{g})}$ is exactly the same as that for the case $\mathfrak{g}=\mathfrak{sl}_2$ in \cite[Proposition 2.5]{MY}, provided that every Weyl module $W^k_{M_\lambda}$ has finite length, where $M_\lambda$ is the finite-dimensional simple $\mathfrak{g}$-module of highest weight $\lambda$. This is proved in \cite[Proposition 2.14]{KL1}, where the proof does not require $\mathfrak{g}$ to be simply laced, only that $k+h^\vee\notin\QQ_{\geq 0}$. So $KL^k(\mathfrak{g})$ is the same as the category of $C_1$-cofinite $V^k(\mathfrak{g})$-modules.

Further, $V^k(\mathfrak{g})$ is a simple vertex operator algebra because $W^k_{M_0}$ is a simple $\widehat{\mathfrak{g}}$-module by \cite[Proposition 2.12]{KL1}, which again does not require $\mathfrak{g}$ to be simply laced, only that $k+h^\vee\notin\QQ$. One can also easily show that $(L^k_\lambda)'\cong L^k_{\lambda^*}$ for any dominant integral weight $\lambda$, where $M_{\lambda^*}\cong M_\lambda^*$. It follows that $V^k(\mathfrak{g})$ is self-contragredient and that $KL^k(\mathfrak{g})$ is closed under contragredients. Now $KL^k(\mathfrak{g})$ is a braided ribbon category by Theorem \ref{thm:C1-cofinite-rigidity}.
\end{proof}

\begin{rem}
For $k+h^\vee\in\QQ_{>0}$, it is best to define $KL^k(\mathfrak{g})$ as the category of finitely-generated grading-restricted generalized $V^k(\mathfrak{g})$-modules. This way, $KL^k(\mathfrak{g})$ still equals the category of $C_1$-cofinite grading-restricted generalized $V^k(\mathfrak{g})$-modules by the argument of \cite[Proposition 2.5]{MY}, so $KL^k(\mathfrak{g})$ is still a braided monoidal category by \cite{Hu}. However, $V^k(\mathfrak{g})$ is no longer generally simple or self-contragredient, so $KL^k(\mathfrak{g})$ is not usually rigid (see \cite{MY} for a detailed analysis of the case $\mathfrak{g}=\mathfrak{sl}_2$). To get a rigid tensor category when $k+h^\vee\in\QQ_{>0}$, one must instead take the smaller Kazhdan-Lusztig category $KL_k(\mathfrak{g})$ of finitely-generated (equivalently, $C_1$-cofinite) grading-restricted generalized modules for the simple affine vertex operator algebra $L_k(\mathfrak{g})$, which is the unique simple quotient of $V^k(\mathfrak{g})$ and is a vertex operator algebra structure on $L^k_0$. Then $L_k(\mathfrak{g})$ is simple and self-contragredient, but it is not clear that $KL_k(\mathfrak{g})$ is closed under contragredients. We will prove this, and show further that $KL_k(\mathfrak{g})$ is a semisimple braided ribbon category for $k+h^\vee\in\QQ_{>0}$, in \cite{McR3}.
\end{rem}

\end{document}